\documentclass[letterpaper,12pt]{article}
\usepackage{amsfonts,latexsym,amssymb,amsmath,epsfig,rotating}
\newlength{\jmr}
\newlength{\bjorn}
\newtheorem{conj}{Conjecture}
\renewcommand{\mod}{\mathrm{mod}}
\newtheorem{agp}{The AGP Theorem}

\newtheorem{yu}{Yu's Theorem.}

\newtheorem{hensel}{Hensel's Lemma}

\newcommand{\sat}{\mathtt{3CNFSAT}}

\renewenvironment{bmatrix}{\left[\begin{array}{*{20}{c}}}{\end{array}\right]}

\newtheorem{dfn}{Definition}[section]
\newtheorem{algor}[dfn]{Algorithm}
 
\newtheorem{rem}[dfn]{Remark}

\newtheorem{prop}[dfn]{Proposition}
\newtheorem{thm}[dfn]{Theorem} 
\newtheorem{lemma}[dfn]{Lemma}
\newtheorem{cor}[dfn]{Corollary}
\newtheorem{ex}[dfn]{Example}

\newcommand{\thth}{^{\text{\underline{th}}}}
\newcommand{\gln}{{\mathbf{GL}_n}}
\newcommand{\rd}{^{\text{\underline{rd}}}}

\newcommand{\nd}{^{\text{\underline{nd}}}}

\newcommand{\ord}{{\mathrm{ord}}}

\newcommand{\np}{{\mathbf{NP}}}

\newcommand{\feas}{{\text{{\tt FEAS}}}}

\newcommand{\fqp}{{\feas_{\Q_\mathrm{primes}}}}

\newcommand{\zpp}{{\mathbf{ZPP}}}

\newcommand{\pp}{\mathbf{P}}

\newcommand{\nc}{\mathbf{NC}}
\newcommand{\expt}{\mathbf{EXPTIME}}
\newcommand{\pspa}{\mathbf{PSPACE}}

\newcommand{\eps}{\varepsilon}
\newcommand{\Pro}{{\mathbb{P}}}
\newcommand{\F}{\mathbb{F}}
\newcommand{\Q}{\mathbb{Q}}
\newcommand{\R}{\mathbb{R}}
\newcommand{\C}{\mathbb{C}}
\newcommand{\N}{\mathbb{N}}
\newcommand{\Z}{\mathbb{Z}}
\newcommand{\bO}{\mathbf{O}}

\newcommand{\Zn}{\Z^n}

\newcommand{\cA}{{\mathcal{A}}}

\newcommand{\cD}{{\mathcal{D}}}
\newcommand{\cE}{{\mathcal{E}}}

\newcommand{\cF}{{\mathcal{F}}}

\newcommand{\cI}{{\mathcal{I}}}

\newcommand{\cL}{{\mathcal{L}}}
\newcommand{\cM}{{\mathcal{M}}}

\newcommand{\cP}{{\mathcal{P}}}
\newcommand{\cQ}{{\mathcal{Q}}}
\newcommand{\cR}{{\mathcal{R}}}
\newcommand{\cS}{{\mathcal{S}}}
\newcommand{\cT}{{\mathcal{T}}}

\newcommand{\cW}{{\mathcal{W}}}

\newcommand{\qed}{$\blacksquare$}
\newcommand{\dia}{$\diamond$}

\newcommand{\newt}{\mathrm{Newt}}

\newcommand{\size}{\mathrm{size}}

\newcommand{\res}{\mathrm{Res}}
\newcommand{\supp}{\mathrm{Supp}}

\newlength{\hwl}
\begin{document}

\title{\mbox{}\\
\vspace{-1in} 
Faster $p$-adic Feasibility for Certain 
Multivariate Sparse Polynomials} 
\author{ 
\scalebox{.8}[1]{Mart\'\i{}n Avenda\~{n}o}\thanks{TAMU 3368, Math Dept., 
College Station, TX \ 77843-3368, USA. {\tt mavendar@yahoo.com.ar} , {\tt 
rojas@math.tamu.edu} , 
{\tt korben@rusek.org} .  Partially supported by NSF MCS grant DMS-0915245.  
J.M.R.\ and K.R.\ also partially supported by Sandia National
Labs and DOE ASCR grant DE-SC0002505. Sandia is a multiprogram
laboratory operated by Sandia Corp., a Lockheed Martin Company, for
the US DOE 
under Contract DE-AC04-94AL85000.} \and 
\scalebox{.8}[1]{Ashraf Ibrahim}\thanks{TAMU 3141,  
Aerospace Engineering Dept., College Station, TX \ 77843-3141, USA,   
{\tt ibrahim@aero.tamu.edu}} \and 
\scalebox{.8}[1]{J. Maurice Rojas}$^*$ \and 
\scalebox{.8}[1]{Korben Rusek}$^*$ } 

\date{\today} 

\maketitle

\begin{abstract} 
\scalebox{.96}[1]{We present algorithms revealing new families of polynomials 
allowing sub-exponential}\linebreak 
detection of $p$-adic rational roots, relative to the sparse 
encoding. For instance, we show that the case of honest $n$-variate 
$(n+1)$-nomials is doable in $\np$ and, for $p$ exceeding 
the Newton polytope volume and not dividing any coefficient, in constant time. 
Furthermore, using the theory of linear forms in 
$p$-adic logarithms, we prove that the case of trinomials in 
one variable can be done in $\np$. The best previous complexity bounds for 
these problems were $\expt$ or worse. Finally, we prove that detecting 
$p$-adic rational roots for sparse polynomials in one variable is $\np$-hard 
with respect to randomized reductions. The last proof makes use of 
an efficient construction of primes in certain arithmetic progressions.  
The smallest $n$ where detecting $p$-adic rational roots for $n$-variate 
sparse polynomials is $\np$-hard appears to have been unknown. 
\end{abstract} 

\section{Introduction}  
Paralleling earlier results over the real numbers \cite{brs}, we study the 
complexity of detecting $p$-adic rational roots for sparse polynomials. 
We find complexity lower bounds over $\Q_p$ hitherto 
unattainable over $\R$, as well as new algorithms over $\Q_p$ 
with complexity close to that of recent algorithms over $\R$ (see 
Theorem \ref{thm:qp} below). 

More precisely, for any commutative ring $R$ with multiplicative identity, 
we let $\feas_R$ --- the {\em $R$-feasibility} 
{\em problem} (a.k.a.\ Hilbert's Tenth Problem over $R$ \cite{h10}) 
--- denote the problem of deciding whether an input 
polynomial system $F\!\in\!\bigcup_{k,n\in\N} (\Z[x_1,\ldots,x_n])^k$ has 
a root in $R^n$.  Observe that $\feas_\R$, $\feas_\Q$, and 
$\{\feas_{\F_q}\}_{q \text{ a prime power}}$ 
are central problems respectively in algorithmic real 
algebraic geometry, algorithmic number theory, and 
cryptography. 

Algorithmic results over the $p$-adics are useful in 
many computational areas: polynomial-time factoring algorithms over $\Q[x_1]$ 
\cite{lll}, computational complexity \cite{antsv}, 
studying prime ideals in 
number fields \cite[Ch.\ 4 \& 6]{cohenant}, elliptic 
curve cryptography \cite{lauder}, and the computation of zeta functions 
\cite{denefver,lauderwan,chambert}. Also, much work has gone into using
$p$-adic methods to algorithmically detect
rational points on algebraic plane curves via variations of
the {\em Hasse Principle}\footnote{
If $F(x_1,\ldots,x_n)\!=\!0$ is any polynomial equation and
$Z_K$ is its zero set in $K^n$, then the Hasse Principle is the
assumption that [$Z_\C$ smooth, $Z_\R\!\neq\!\emptyset$, and
$Z_{\Q_p}\!\neq\!\emptyset$ for all primes $p$] implies
$Z_\Q\!\neq\!\emptyset$ as well. The Hasse Principle is a theorem when
$Z_\C$ is a quadric hypersurface or a curve of genus zero,
but fails in subtle ways already for curves of genus one (see, e.g.,
\cite{bjornhasse1}). } (see, e.g., \cite{colliot,bjornbm}). 
However, our knowledge of the complexity of deciding the existence of 
solutions for {\em sparse} polynomial equations over $\Q_p$ is surprisingly 
coarse: good bounds for the number of solutions over $\Q_p$ in 
one variable weren't even known until the late 1990s \cite{lenstra2}. 
\begin{dfn}
\label{dfn:qp}
Let $\fqp$ denote the problem of deciding, for an input Laurent polynomial
system $F$ $\in\!\bigcup_{k,n\in\N}\left(\Z\!\left[x^{\pm 1}_1,\ldots,
x^{\pm 1}_n \right]\right)^k$ {\em and}
an input prime $p$, whether $F$ has a root in $\Q^n_p$. 
Also let $\Pro\!\subset\!\N$ denote the set of primes, $p\!\in\!\Pro$, and, 
when $\cI$ is a family of such pairs $(F,p)$, we let $\fqp(\cI)$ denote the 
restriction of $\fqp$ to inputs in $\cI$. 

When $a_j\!\in\!\Zn$, the notations $a_j\!=\!(a_{1,j},\ldots,a_{n,j})$,
$x^{a_j}\!=\!x^{a_{1,j}}_1\cdots x^{a_{n,j}}_n$, and $x\!=\!(x_1,\ldots,x_n)$
will be understood. Also, when $f(x)\!:=\!\sum^m_{j=1} c_ix^{a_j}$
with $c_j\!\in\!\Z\setminus\{0\}$ for all $j$,
and the $a_j\!\in\!\Zn$ are pair-wise distinct, we call $f$ an 
{\em $n$-variate $m$-nomial}, and we define
$\supp(f)\!:=\!\{a_1,\ldots,a_m\}$ to be the {\em support} of $f$. 
We also define $\newt(f)$ --- the (standard) Newton polytope of 
$f$ --- to be the convex hull of\footnote{i.e., smallest convex set
containing...} $\supp(f)$ and let $V_f$ denote its $n$-dimensional 
volume, normalized so that $[0,1]^n$ has volume $1$. 

Let $\size(f)\!:=\!\sum^m_{i=1}
\log_2\left[(2+|c_i|)(2+|a_{1,i}|)\cdots
(2+|a_{n,i}|)\right]$ and $\size(F)\!:=\!\sum^k_{i=1}\size(f_i)$.\linebreak  
\scalebox{.91}[1]{The underlying input sizes for $\fqp$ and $\fqp(\cI)$ shall 
then be $\size_p(F)\!:=\!\size(F)+\log p$,}\linebreak   
and we use $\size(F)$ as the input size for $\feas_{\Q_p}$ for any 
prime $p$. Finally, we let $\cF_{n,m}$ denote the set of all $n$-variate
$m$-nomials and, for any $m\!\geq\!n+1$, we let
$\cF^*_{n,m}\!\subseteq\!\cF_{n,m}$
denote the subset consisting of those $f$ with $V_f\!>\!0$ 
We call any $f\!\in\!\cF^*_{n,m}$ an {\em 
honest $n$-variate $m$-nomial} (or {\em honestly $n$-variate}). \dia
\end{dfn}

\noindent 
As an example, it is clear that upon 
substituting $y_1\!:=\!x^2_1x_2x^7_3x^3_4$, 
the dishonestly $4$-variate trinomial
$-1+7x^2_1x_2x^7_3x^3_4-43x ^{198}_1x^{99}_2x^{693}_3x^{297}_4$
(with support contained in a line segment) has a root 
in $(\Q^*_p)^4$ iff the {\em honest uni}variate trinomial
\mbox{$-1+7y_1-43y^{99}_1$} has a root in 
$\Q^*_p$. Via the use of Hermite Normal Form (as in Section \ref{sec:binom} 
below), it is then easy to see that there is no loss of generality in 
restricting to $\cF^*_{n,n+k}$ (with $k\!\geq\!1$) when studying the 
algorithmic complexity of sparse polynomials. 
Note also that the degree, $\deg f$, of a polynomial $f$ can sometimes be
exponential in $\size(f)$ for certain families of $f$, e.g.,
$d\!\geq\!2^{\size\left(1+5x^{126}_1+x^d_1\right)-16}$.

While there are now randomized algorithms for factoring $f\!\in\!\Z[x_1]$ 
over $\Q_p[x_1]$ with expected complexity polynomial in $\deg(f)+\size_p(f)$ 
\cite{cantorqp} (see also \cite{chistov}), no such algorithms are known   
to have complexity polynomial in $\size_p(f)$ alone. Our main theorem below 
shows that such algorithms are hard to derive because finding just the linear 
factors is already essentially equivalent to the 
$\pp\text{\scalebox{1}[.85]{$\stackrel{?}{=}$}}\np$ problem. Nevertheless, 
we obtain fast new algorithms for interesting sub-cases of 
$\fqp\!\left(\bigcup_{n\in\N}\Z[x_1,\ldots,x_n]\right)\times\Pro)$.  
\begin{thm}
\label{thm:qp} \mbox{}\\
\scalebox{.88}[1]{{\bf 0.} $\fqp(\cF_{1,m}\times \Pro)\!\in\!\pp$ for 
$m\!\in\!\{0,1,2\}$. \ \ \ \  
{\bf 1.} For any fixed prime $p$ we have 
$\feas_{\Q_p}(\cF_{1,3})\!\in\!\np$.}\linebreak
{\bf 2.} There is a countable union of algebraic hypersurfaces
$\cE\; \subsetneqq\; \Z[x_1]\times\Pro$, with natural density\linebreak 
\mbox{}\hspace{.6cm}$0$, such that
$\fqp((\Z[x_1]\times \Pro)\setminus \cE)\!\in\!\np$.\\
{\bf 3.} (a) $\fqp\!\left(\left(\bigcup_{n\in\N} \cF^*_{n,n+1}\right)\times 
\Pro\right)\!\in\!\np$.\\ 
\mbox{}\hspace{.5cm}(b) \scalebox{.87}[1]{Letting
$\cQ\!:=\!\{c_0+c_1x^2_1+\cdots+c_nx^2_n\; | \; n\!\in\!\N; 
\ c_0,\ldots,c_n \!\in\!\Z\setminus\{0\}\}\times \Pro$, we have 
$\fqp(\cQ)\!\in\!\pp$.}\\
\mbox{}\hspace{.5cm}(c) Letting $\cW\!\subset\!\left(\bigcup_{n\in 
\N} \cF^*_{n,n+1}\right)\times\Pro$ denote the subset consisting of 
those $(f,p)$ with $n\!\geq\!2$,\linebreak 
\mbox{}\hspace{1.2cm}$p\!\geq\!(n!V_F)^{2/(n-1)}$, and $p$ not 
dividing $n!V_F$ or any coefficient of $f$, we have\linebreak 
\mbox{}\hspace{1.2cm}that $f$ always has a root in $\Q^n_p$ for any 
$(f,p)\!\in\!\cW$, i.e., $\fqp(\cW)$ is doable in\linebreak 
\mbox{}\hspace{1.2cm}constant time.\\
{\bf 4.} If $\fqp(\Z[x]\times \Pro)\!\in\!\zpp$ then $\np\!\subseteq\!\zpp$.\\
{\bf 5.} If the Wagstaff Conjecture is true, then
$\fqp(\Z[x]\times \Pro)\!\in\!\pp \Longrightarrow \pp\!=\!\np$,
i.e., we\linebreak 
\mbox{}\hspace{.6cm}can strengthen Assertion (4) above. 
\end{thm} 

\noindent 
The aforementioned complexity classes, are reviewed 
briefly in Section \ref{sec:back} (see also \cite{papa} for an excellent 
textbook treatment). The {\em Wagstaff Conjecture}, dating back to 1979
(see, e.g., \cite[Conj.\ 8.5.10, pg.\ 224]{bs}), is the 
assertion that the least prime congruent to $k$ mod $N$ is 
$O(\varphi(N)\log^2 N)$,
where $\varphi(N)$ is the number of integers in $\{1,\ldots,N\}$ relatively
prime to $N$. This conjectural bound is (unfortunately) much stronger than the 
known implications of the Generalized Riemann Hypothesis. 

Let us now briefly highlight what is new in our main theorem, and how 
the real case compares.\footnote{A weaker version of Theorem 
\ref{thm:qp}, without Assertions (1) and (3), 
appeared recently in an extended abstract \cite{padic1}. } 
First, one can in fact prove $\feas_\R\!\left(\bigcup_{n\in\N} \cF^*_{n,n+1}
\right)\!\in\!\nc^1$ (i.e., a much stronger real analogue of Assertion (3)) 
via some elementary tricks involving monomial changes of variables 
\cite[Thm.\ 1.3]{brs}. 
Unfortunately, these tricks are obstructed over $\Q_p$ 
(see Example \ref{ex:martin} below), thus making Assertion (3) harder to 
prove. As evinced by Parts (b) and (c) of Assertion (3), algorithms 
for $\fqp\!\left(\left(\bigcup_{n\in\N} \cF^*_{n,n+1}\right)\times 
\Pro\right)$ clearly complement classical results on quadratic forms 
(see, e.g., \cite[Ch.\ IV]{serre}) and the Weil Conjectures (see, e.g., 
\cite{weil,freitag}). More to the point, 
the best previous complexity upper bound for 
$\fqp\!\left(\left(\bigcup_{n\in\N} \cF^*_{n,n+1}\right)\times\Pro\right)$ 
appears to be quadruply exponential, via an extension 
of Hensel's Lemma by Birch and McCann \cite{birchmccann}.  

While the real analogue of Assertion (0) is not hard to prove, 
$\feas_\R(\cF_{1,3})\!\in\!\pp$ (a stronger real analogue for Assertion (1)) 
was proved only recently \cite[Thm.\ 1.3]{brs} using linear forms in 
logarithms \cite{nesterenko}. It is thus worth noting 
that the proof of Assertion (1) (in Section \ref{sec:lin}) uses linear forms 
in $p$-adic logarithms \cite{yu} at a critical juncture, and suggests an 
approach to a significant speed-up. 
\addtocounter{footnote}{1}
\begin{cor}
\label{cor:faster} 
Suppose that for all $p\!\in\!\Pro$ and $\ell\!\geq\!1$, 
$\feas_{\Z/p^\ell\Z}(\cF_{1,3})$ admits a 
(deterministic) algorithm$^3$ with 
complexity $(p+\ell+\size(f))^{O(1)}$.\footnotetext{All algorithms 
discussed here are based on Turing machines \cite{papa}.} 
Then for any fixed prime $p$, $\feas_{\Q_p}(\cF_{1,3})\!\in\!\pp$. 
\end{cor} 

\noindent 
The truth of the hypothesis to our corollary above appears to be 
an open question. (Note that brute-force search easily leads to 
an algorithm of complexity $p^\ell \size(f)^{O(1)}$, so the main issue 
here is the dependence on $\ell$.) Paraphrased in our notation, 
Erich Kaltofen asked in 2003 whether $\feas_{\Z/p\Z}(\cF_{1,3})$ admits a 
(deterministic) algorithm with complexity\linebreak  
$(\log(p)+\size(f))^{O(1)}$ \cite{kaltofen}.\footnote{David A.\ Cox also 
independently asked Rojas the same question in august of 2004.}   

The best previous complexity upper bound for $\fqp(\Z[x_1]\times \Pro)$ 
relative to the sparse input size appears to have been $\expt$ 
\cite{mw}. In particular,\linebreak 
$\fqp(\cF_{1,4}\times\Pro)\text{\scalebox{1}[.8]{$\stackrel{?}{\in}$}}
\np$ and $\feas_\R(\cF_{1,4})\text{\scalebox{1}[.8]{$\stackrel{?}{\in}$}}\np$
are still open questions \cite[Sec.\ 1.2]{brs}. High probability 
speed-ups over $\R$ paralleling Assertion (2) are also unknown at this time.  
For clarity, here is an example illustrating the zero-density exception in 
Assertion (2).  
\begin{ex}
\label{ex:lots}
\scalebox{.92}[1]{Let $T$ denote the family of pairs
$(f,p)\!\in\!\Z[x_1]\times \Pro$ with $f(x_1)\!=\!a+bx^{11}_1+cx^{17}_1
+x^{31}_1$}\linebreak 
and let $T^*\!:=\!T\setminus \cE$. Then there is
a sparse $61\times 61$ structured matrix $\cS$ (cf.\ Lemma \ref{lemma:syl}
in Section \ref{sec:disc} below), whose entries lie in
$\{0,1,31,a,b,11b,c,17c\}$,
such that $(f,p)\!\in\!T^* \Longleftrightarrow$\linebreak 
\mbox{$p\not|\det \cS$}.
So by Theorem \ref{thm:qp}, $\fqp(T^*)\!\in\!\np$, and Corollary \ref{cor:lots}
in Section \ref{sec:disc} below tells us that for large coefficients,
$T^*$ occupies almost all of $T$. In particular,
letting $T(H)$ (resp.\ $T^*(H)$) denote those
pairs $(f,p)$ in $T$ (resp.\ $T^*$) with
$|a|,|b|,|c|,p\!\leq\!H$,
we obtain\\
\mbox{}\hfill $\frac{\#T^*(H)}{\#T(H)}\!\geq\!\left(1-\frac{244}{2H+1}\right)
\left(1-\frac{1+61\log(4H)\log H}{H}\right)$. \hfill \mbox{}\\
In particular, one can check via {\tt Maple} that\\
\mbox{}\hfill
$(-973+21x^{11}_1-2x^{17}_1 +x^{31}_1,p)\!\in\!T^*$\hfill\mbox{}\\
for all but $352$ primes $p$. \dia
\end{ex}

As for lower bounds, the least $n$ making 
$\fqp(\Z[x_1,\ldots,x_n]\times \Pro)$ $\np$-hard appears to have been 
unknown. Assertions (4) and (5) thus come close to settling this problem. 
In particular, while is not hard to show that the full problem $\fqp$ is 
$\np$-hard, the proofs of Assertions (4) and (5) make essential use of a 
deep result of Alford, Granville, and Pomerance 
\cite{carmichael} on primes in random arithmetic progressions.  
We detail this connection below. 

\subsection{Related Work, a Topological Observation, Weil's Conjecture, 
and Primes in Arithmetic Progression} 
\label{sub:rel}
Let us first recall that Emil Artin conjectured around 1935 that, for any 
prime $p$,  homogeneous polynomials of degree $d$ in $n\!>\!d^2$ variables 
always have non-trivial roots in $\Q^n_p$ \cite{artin}. (The 
polynomials $x^2_1+\cdots+x^2_n$ show that Artin's conjecture 
is resoundingly false over the real numbers.) 
Artin's conjecture was already known to be true for $d\!=\!2$ 
\cite{hasse} and, in 1952, the $d\!=\!3$ case was proved by Lewis \cite{lewis}. 
However, in 1966, Terjanian 
disproved the conjecture via an example with $(p,d,n)\!=\!(2,4,18)$. 

The Ax-Kochen Theorem from 1965 provided a valid correction of 
Artin's conjecture: for any $d$, there is a constant $p_d$ such that for all 
primes $p\!>\!p_d$, any homogeneous degree $d$ polynomial 
in $n\!>\!d^2$ variables has a $p$-adic rational root \cite{axkochen,
heathbrown}. The hard cases of $\fqp$ then appear to consist of
high degree polynomials with few variables and $p$ small.  

It is interesting to observe that while it is easier for a 
polynomial in many variables to have roots over $\Q_p$ than over $\R$, 
deciding the existence of roots appears to be much harder over $\Q_p$ than 
over $\R$. In particular, while Tarski showed in 1939  
that $\feas_\R$ is decidable \cite{tarski}, $\feas_{\Q_p}$ wasn't shown to be 
decidable until work of Cohen in the 1960s \cite{cohen}. Now, the 
best general complexity upper bounds appear to be $\pspa$ for $\feas_\R$ 
\cite{pspace} and {\em quadruply} exponential for $\feas_{\Q_p}$ 
\cite{birchmccann,greenberg}.  

While the univariate problems $\feas_\R(\cF_{1,2})$ and 
$\fqp(\cF_{1,2})$ are now both known to be in $\pp$, their natural multivariate 
extensions $\feas_\R\!\left(\bigcup_{n\in \N}\cF^*_{n,n+1}\right)$ 
and $\feas_{\Q_p}\!\left(\bigcup_{n\in \N}\cF^*_{n,n+1}\right)$ 
already carry nuances distinguishing the real and $p$-adic settings: 
topological differences between the 
real and $p$-adic zero sets of polynomials in $\cF^*_{n,n+1}$ 
force the underlying feasibility algorithms to differ. 
Concretely, positive zero sets for polynomials in $\cF^*_{n,n+1}$ 
are always either empty or non-compact. This in turn allows one to solve 
$\feas_\R\!\left(\bigcup_{n\in \N}\cF^*_{n,n+1}\right)$ by simply 
checking signs of coefficients, independent of the exponents 
\cite[Thm.\ 1.3]{brs}. On the other hand, solving 
$\feas_{\Q_p}\!\left(\bigcup_{n\in \N}\cF^*_{n,n+1}\right)$ depends
critically on the exponents (see Corollary \ref{cor:binom} of 
Section \ref{sec:binom}), and the underlying hypersurfaces in $\Q^n_p$ can 
sometimes be a single isolated point.  
\begin{ex} 
\label{ex:martin} 
Consider $f(x_1,x_2)\!:=\!1+2x^2_1-3x^2_2$. 
Then it is easy to see that $(1,1)$ is the unique root of $f$ in 
$\F^2_7$. Via Hensel's Lemma (see Section \ref{sec:back} below), 
the root $(1,1)\!\in\!\F^2_7$ can then 
be lifted to a unique root of $f$ in $\Q^2_7$. In particular, by checking 
valuations, any root of $f$ in $\Q^2_7$ must be the lift of some root of $f$ 
in $\F^2_7$, and thus $(1,1)$ is the only root of $f$ in $\Q^2_7$. \dia 
\end{ex}  

Our last example illustrated the importance of finite fields in studying 
$p$-adic rational roots. Deligne's Theorem on zeta functions over 
finite fields (n\'ee the Weil Conjectures) is the definitive statement 
on the connection between point counts over finite fields and complex 
geometry. The central result that originally motivated the  
Weil Conjectures will also prove useful in our study of $\fqp$.  
\begin{thm} 
\label{thm:weil} \cite[Pg.\ 502]{weil}
Let $p$ be any prime, $d_1,\ldots,d_n\!\in\!\N$, and let 
$c_0,\ldots,c_n$ be integers not divisible by $p$. 
Then, defining $f(x)\!:=\!c_0+c_1x^{d_1}_1+\cdots+c_nx^{d_n}$, the 
number, $N$, of roots of $f$ in $\F^n_p$ satisfies 
$|N-p^{n-1}|\leq\left(\prod^n_{i=1} (\gcd(d_i,p-1)-1)\right)p^{(n-1)/2}$. \qed 
\end{thm} 

Finally, it is worth noting that our $\np$-hardness proof requires 
the efficient construction of primes in certain arithmetic progressions. 
The following result, inspired by earlier work of von zur Gathen, 
Karpinski, and Shparlinski, may be of independent interest. 
\begin{thm} 
\label{thm:primes} 
For any $\delta\!>\!0$,
$\eps\!\in\!(0,1/2)$, and $n\!\in\!\N$, we can find ---
within\\
\mbox{}\hfill $O\!\left((n/\eps)^{\frac{3}{2}+\delta} + 
(n\log(n)+\log(1/\eps))^{7+\delta}\right)$ \hfill \mbox{}\\
randomized bit operations --- a sequence $P\!=\!(p_i)^n_{i=1}$ of consecutive 
primes and $c\!\in\!\N$ such that $p\!:=\!1+c\prod^n_{i=1}p_i$
satisfies $\log p = O(n\log(n)+\log(1/\eps))$
and, with probability $\geq\!1-\eps$, $p$ is prime.  
\end{thm} 

\subsection{Future Directions} 
Since $\np$-hardness is easier to prove for detecting roots 
of univariate polynomials over $\Q_p$ than over $\R$, we anticipate 
that a similar phenomenon occurs for multivariate polynomials. 
\begin{conj} 
For any fixed prime $p$ we have that 
$\feas_{\Q_p}\!\left(\bigcup_{n\in\N} \cF^*_{n,n+1}\right)$ 
is $\np$-hard. 
\end{conj} 

\noindent 
It is already known that 
$\feas_\R\!\text{\scalebox{1}[.8]{$\left(\bigcup
\limits_{n\in\N \ , \ 0<\eps'\leq\eps}
\cF^*_{n,n+n^{\eps'}}\right)$}}$ is $\np$-hard for any $\eps\!>\!0$ 
\cite[Thm.\ 1.3]{brs}. In particular, it is likely one can modify the proof 
of the latter statement to at least 
prove that $\feas_{\Q_p}\!\text{\scalebox{1}[.8]{$\left(\bigcup
\limits_{n\in\N \ , \ 0<\eps'\leq\eps}
\cF^*_{n,n+n^{\eps'}}\right)$}}$ is $\np$-hard for any fixed prime $p$.  

Further speed-ups for detecting $p$-adic rational roots of 
$n$-variate $(n+1)$-nomials appear to hinge on a better understanding 
of the analogous problem over certain finite rings. In particular, the 
truth of the following conjecture would imply 
$\fqp\!\left(\cF^*_{n,n+1}\right)\!\in\!\pp$ for any fixed $n$. 
\begin{conj} 
\label{conj:simplex} 
Suppose $\ell,n\!\in\!\N$ and $p\!\in\!\Pro$. Then 
$\feas_{\Z/p^\ell\Z}(\cF^*_{n,n+1})$  admits a (deterministic) algorithm 
with complexity $(\log(p)+\ell+\size(f))^{O(n)}$.
\end{conj} 

\noindent 
Note that brute-force search easily attains a complexity bound of 
$p^{\ell n}\size(f)^{O(1)}$ so the key difficulty is the dependence on 
$p^\ell$. 

Finally, it is worth noting that 
$\feas_\R(\cF^*_{n,n+2})\!\in\!\pp$ for any fixed $n\!\in\!\N$ 
\cite[Thm.\ 1.3]{brs}. In fact, the proof there inspired our proof of 
Assertion (1) of Theorem \ref{thm:qp}, so it would be most interesting to 
extend our techniques to the multivariate case. 
\begin{conj} 
For any fixed $n\!\in\!\N$ and $p\!\in\!\Pro$ we have 
$\feas_{\Q_p}(\cF^*_{n,n+2})\!\in\!\np$. 
\end{conj} 

We review some general background in Section \ref{sec:back} before proving our 
main results.  Some of the results we'll need will appear just before their use
in the proofs of Assertions (0) and (3) in Section \ref{sec:binom}, 
the proof of Assertion (2) in Section \ref{sec:disc}, the 
proof of Assertion (1) in Section\linebreak 
\scalebox{.95}[1]{\ref{sec:lin}, the proof of 
Theorem \ref{thm:primes} in Section \ref{sub:agp}, and the proofs of 
Assertions (4) and (5) in Section \ref{sec:primes}.} 

\section{Complexity Classes and $p$-adic Basics} 
\label{sec:back} 
Let us first recall briefly the following complexity classes
(see also \cite{papa} for an excellent textbook treatment):
\begin{itemize}
\item[$\nc^1$]{The family of functions computable by Boolean
circuits with size polynomial\footnote{Note that the underlying polynomial
depends only on the problem in question (e.g., matrix inversion, shortest path
finding, primality detection) and not the particular instance of the
problem.} in the input size and depth
$O(\log^i \text{{\tt InputSize}})$. }
\item[$\pp$]{ The family of decision problems which can be done within time
polynomial in the input size. }
\item[$\zpp$]{ The family of decision problems admitting a
randomized polynomial-time algorithm giving a correct answer,
or a report of failure, the latter occuring with probability
$\leq\!\frac{1}{2}$. }
\item[$\np$]{ The family of decision problems where a ``{\tt Yes}'' answer can
be {\em certified} within time polynomial in the input size.} 
\item[$\pspa$]{ The family of decision problems solvable within time
polynomial in the input size, provided a number of processors exponential in
the input size is allowed. }
\item[$\expt$]{ The family of decision problems solvable
within time exponential in the input size.}
\end{itemize}

The following containments are standard:\\
\mbox{}\hfill
$\nc^1\subseteq\pp\subseteq\zpp\subseteq\np\subseteq\pspa\subseteq
\expt$.\hfill\mbox{}\\
The properness of each adjacent inclusion above (and even 
the properness of $\pp\!\subseteq\!\pspa$)  
is a major open problem \cite{papa}.

Recall that for any ring $R$, we denote its unit group by $R^*$.
For any prime $p$ and $x\!\in\!\Z$, recall that the
{\em $p$-adic valuation}, $\ord_p x$, is the greatest $k$ such that $p^k|x$.
We can extend $\ord_p(\cdot)$ to $\Q$ by
$\ord_p\frac{a}{b}\!:=\!\ord_p(a)-\ord_p(b)$ for any
$a,b\!\in\!\Z$; and we let $|x|_p\!:=\!p^{-\ord_p x}$ denote the {\bf
$p$-adic norm}. The norm $|\cdot|_p$ defines a natural metric satisfying
the ultrametric inequality and $\Q_p$ is, to put it tersely, the
completion of $\Q$ with respect to this metric. This metric,
along with $\ord_p(\cdot)$, extends naturally to the field of {\em $p$-adic
complex numbers} $\C_p$, which is the metric completion of the algebraic
closure of $\Q_p$ \cite[Ch.\ 3]{robert}.

\scalebox{.94}[1]{It will be useful to recall some classical invariants for 
treating quadratic polynomials over $\Q_p$.} 
\begin{dfn}
\label{dfn:serre}  
\cite[Ch.\ I--IV, pp.\ 3--39]{serre} 
For any prime $p$ and $a\!\in\!\Z$ we define the\linebreak  
{\em Legendre symbol}, $\left(\frac{a}{p}\right)$, to be 
$+1$ or $-1$ according as $a$ has a square root mod $p$ or not. 
Also, for any $b\!\in\!\Z$, we let the {\em ($p$-adic) Hilbert symbol}, 
$(a,b)_p$, be $+1$ or $-1$ according as 
$ax^2+by^2\!=\!z^2$ has a solution in $\Pro^2_{\Q_p}$ or not. 
Finally, for any 
$f(x)\!=\!c_0+c_1x^2_1+\cdots+c_nx^2_n\!\in\!\Z[x_1,\ldots,x_n]$,  
we define $d_f\!:=\!\prod^n_{i=1}c_i$ and $\eps_f\!:=\!\prod_{1\leq i<j\leq 
n} (c_i,c_j)_p$. \dia 
\end{dfn}
\begin{thm} 
\label{thm:quad} 
\cite[Thm.\ 1, pg.\ 20 \& Cor., pp.\ 37]{serre}  
Following the notation of Definition \ref{dfn:serre}, let $j\!:=\!\ord_p a$ 
and $k\!:=\!\ord_p b$. Then the Hilbert symbol $(a,b)_p$ is exactly \\
\mbox{}\hspace{1cm}(i) $(-1)^{jk(\frac{p-1}{2} \mod \ 2)}
\left(\frac{a/p^j}{p}\right)^k  
\left(\frac{b/p^k}{p}\right)^j$, or \\
\mbox{}\hspace{1cm}(ii) $(-1)^{Z(a,b)}$ where 
$Z(a,b):=\left(\frac{a/p^{j}-1}{2}\right) \left(\frac{b/p^{k}-1}{2}\right) 
             + j\left(\frac{(b/p^k)^2-1}{8}\right)
             + k\left(\frac{(a/p^j)^2-1}{8}\right) \ \mod \ 2$,\\ 
according as $p\!\neq\!2$ or $p\!=\!2$. 

Finally, $f$ has a root in $\Q_p$ iff one of the following conditions 
holds:\\
1.\ $n\!=\!1$, $\mu\!:=\!\ord_p(c_0/c_1)$ is even, and 
$\left(\frac{-c_0/(c_1p^\mu)}{p}\right)\!=\!1$. \\ 
2.\ $n\!=\!2$ and $(-c_0,-d_f)_p\!=\!\eps_f$ (viewing $c_0$
and $d_f$ as elements of $\Q_p/(\Q^*_p)^2$). \\
3.\ $n\!=\!3$ and either $c_0\!\neq\!d_f$ or $c_0\!=\!d_f$ and 
$(-1,-d_f)\!=\!\eps_f$ (viewing $c_0$
and $d_f$ as elements of\linebreak
\mbox{}\hspace{.5cm}$\Q_p/(\Q^*_p)^2$). \\
4.\ $n\!\geq\!4$. \qed 
\end{thm} 

A key tool we will use throughout this paper is 
{\em Hensel's Lemma}, suitably extended to multivariate Laurent polynomials. 
\begin{hensel}
Suppose $f\!\in\!\Z_p\!\left[x^{\pm 1}_1,\ldots,x^{\pm 1}_n\right]$ 
and $\zeta_0\!\in\!\Z^n_p$ satisfies 
$\ord_p \frac{\partial f}{\partial x_i}(\zeta_0)\!=\!\ell\!<\!
+\infty$ for some $i\!\in\!\{1,\ldots,n\}$, and   
$f(\zeta_0)\!\equiv\!0 \ (\mod \ p^{2\ell+1})$.  
Then there is a root $\zeta\!\in\!\Z^n_p$ of $f$ with $\zeta\!\equiv\!\zeta_0 \
(\mod \ p^{\ell})$ and $\ord_p \frac{\partial f}{\partial x_i}(\zeta)\!=\!
\ord_p \frac{\partial f}{\partial x_i}(\zeta_0)$. \qed 
\end{hensel}

\noindent 
The special case of polynomials appears as Theorem 1 on the 
bottom of Page 14 of \cite{serre}. (See also \cite{birchmccann}.) 
The proof there extends almost verbatim to Laurent polynomials.  

\section{From Binomials to $(n+1)$-nomials: Proving\\ Assertions (0) and (3)} 
\label{sec:binom} 
Let us first recall the following standard lemma on taking radicals  
in certain finite groups. 
\begin{lemma}
\label{lemma:qp}
(See, e.g., \cite[Thm.\ 5.7.2 \& Thm.\ 5.6.2, pg.\ 109]{bs})
Given any cyclic group $G$, $a\!\in\!G$, and an integer $d$, the following
3 conditions are equivalent:\\
\mbox{}\hspace{1cm}1.\ The equation $x^d\!=\!a$ has a solution.\\
\mbox{}\hspace{1cm}2.\ The order of $a$ divides $\frac{\#G}{\gcd(d,\#G)}$.\\
\mbox{}\hspace{1cm}3.\ $a^{\#G/\gcd(d,\#G)}\!=\!1$.\\
Also, $\F^*_q$ is cyclic for
any prime power $q$, and $(\Z/p^\ell\Z)^*$ is cyclic
for any $(p,\ell)$ with $p$ an odd prime or $\ell\!\leq\!2$.
Finally, for
\scalebox{.95}[1]{$\ell\!\geq\!3$, $(\Z/2^\ell\Z)^*\!=\!
\left\{\pm 1,\pm 5,\pm 5^2,\pm 5^3,\ldots,\pm 5^{2^{\ell-2}-1} \ \mod \ 2^\ell
\right\}$. \qed}
\end{lemma}

A direct consequence of Lemma \ref{lemma:qp} and Hensel's Lemma 
is the following characterization of univariate binomials with 
$p$-adic rational roots. 
\begin{cor} 
\label{cor:binom} 
Suppose $c\!\in\!\Q^*_p$ and $d\!\in\!\Z\setminus\{0\}$. Let 
$k\!:=\!\ord_p c$,  $\ell\!:=\!\ord_p d$, and (if $p\!=\!2$ and 
$d$ is even) $d'\!=\!\left(\frac{d}{2^\ell}\right)^{-1} \ 
(\mod \ 2^{2\ell-1})$. Then the equation $x^d\!=\!c$ has a solution in $\Q_p$ 
iff $d|\ord_p c$ and one of the following two conditions hold:\\
(a) $p$ is odd and $\left(\frac{c}{p^k}
\right)^{p^{\ell}(p-1)}\!=\!1 \ (\mod \ p^{2\ell+1})$.\\ 
(b) \scalebox{.95}[1]{$p\!=\!2$ and either 
(i) $d$ is odd, or (ii) $\left(\frac{c}{p^k}\right)^{d'}\!=\!1 \ (\mod \ 8)$ 
and $\left(\frac{c}{p^k}\right)^{d'2^{\max\{\ell-2,0\}}}\!\!\!\!\!=\!1 
\ (\mod \ 2^{2\ell+1})$.}\\  In particular, these conditions can be 
checked in time polynomial in $\log(d) + \log(p)$ when $\log c \!=\!
(\log(d)+\log(p))^{O(1)}$. Furthermore, when $\ord_p c\!=\!0$, 
$x^d\!=\!c$ has a root in $\Q_p$ iff $x^d\!=\!c$ has a root in 
$(\Z/p^{2\ell+1}\Z)^*$. 
\end{cor} 

\noindent 
{\bf Proof:} 
Replacing $x$ by $1/x$, we can clearly assume $d\!>\!0$. 
Clearly, any $p$-adic root $\zeta$ of $x^d-c$ satisfies
$d\ord_p\zeta\!=\!\ord_p c$. This accounts for the condition
preceding Conditions (a) and (b). 

Replacing $x$ by $p^{\ord_p c/d}x$ (which clearly preserves the existence of
roots in $\Q^*_p$) we can assume further that $\ord_p c\!=\!\ord_p\zeta\!=\!0$. 
Moreover, $\ord_p f'(\zeta)\!=\!\ord_p(d)+(d-1)\ord_p\zeta\!=\!\ord_p d$. So
by Hensel's Lemma, $x^d-c$ has a root in $\Q^*_p$ iff $x^d-c$ has a
root in $(\Z/p^{2\ell+1}\Z)^*$. Lemma \ref{lemma:qp} then immediately
accounts for Condition (a) when $p$ is odd.

Condition (b) then follows routinely: First, one observes that
exponentiating by an odd power is an automorphism of $(\Z/2^{2\ell+1})^*$,
and thus $x^d-c$ has a root in $(\Z/2^{2\ell+1}\Z)^*$ iff
$x^{2^\ell}-c^{d'}$ does. Should $\ell\!=\!0$ then one has a root
regardless of $c$. Otherwise, $c^{d'}$ must be a square for there to be
a root. Since $\ord_p c\!=\!0$, $c$ is odd and  
\cite[Ex.\ 38, pg.\ 192]{bs} tells us that $c^{d'}$
is a square in $(\Z/2^\ell\Z)^*$ iff $c^{d'}\!=\!1 \; (\mod \ 8)$.
Invoking Lemma \ref{lemma:qp} once more on the the cyclic
subgroup $\{1,5^2,5^4,5^6,\ldots,5^{2^{2\ell-1}-2}\}$, it is clear 
that Condition (b) is exactly what we need when $p\!=\!2$. 

To conclude, recall that arithmetic in
$\Z/p^{2\ell+1}\Z$ can be done in time polynomial in $\log(p^\ell)$
\cite[Ch.\ 5]{bs}. Recall also that, in any ring, $x^n$ can be 
computed using just $O(\log n)$ bit operations and  
multiplication of powers of $x$, via recursive squaring 
\cite[Thm.\ 5.4.1, pg.\ 103]{bs}. 
Our conditions are then clearly simple enough to yield the asserted 
time bound. 

The final assertion follows immediately from setting $k\!=\!0$ in 
the conditions we've just derived. \qed

\medskip 
At this point, the proof of Assertion (0) of Theorem \ref{thm:qp} is 
trivial. By combining our last result with a classical integral 
matrix factorization, Assertion (3) then also becomes easy to prove. 
So let us first motivate the connection between 
$n$-variate $(n+1)$-nomials and matrices. 
\begin{prop} 
\label{prop:mat}  
Suppose $K$ is any field, $c_0,\ldots,c_n\!\in\!K$ with $c_i\!\neq\!0$ 
for some $i\!\in\!\{1,\ldots,n\}$, 
$a_1,\ldots,a_n\!\in\!\Zn$ are 
linearly independent vectors, $A$ is the $n\times n$ matrix with 
columns $a_1,\ldots,a_n$, and $f(x)\!:=\!c_0+c_1x^{a_1}+\cdots+c_nx^{a_n}$. 
Then, letting $x\!=\!(x_1,\ldots,x_n)\!\in\!(K^*)^n$ and 
$f_i\!:=\!\frac{\partial f}{\partial x_i}$ for all $i$, we 
have:  
\\
\mbox{}\hfill
$[f_1(x),\ldots,f_n(x)]\!=\![c_1x^{a_1},
\ldots, c_nx^{a_n}] A^T \begin{bmatrix}
x^{-1}_1 &        & \\
           & \ddots & \\
           &        & x^{-1}_n 
\end{bmatrix}$. \hfill \mbox{}\\ 
In particular, all the roots of $f$ in $(K^*)^n$ are 
non-degenerate. 
\end{prop}  

\noindent 
{\bf Proof:} The first assertion is routine. For the second assertion, 
observe that if $\zeta\!\in\!(K^*)^n$ is any root of $f$ then, 
thanks to our first assertion,  
the vector $[f_1(\zeta),\ldots,f_n(\zeta)]$ can not vanish because 
$\det A\!\neq\!0$. \qed 

\begin{dfn}
\label{dfn:smith}
Let $\Z^{n\times n}$ denote the set of $n\times n$ matrices
with all entries integral, and let $\gln(\Z)$ denote the
set of all matrices in $\Z^{n\times n}$ with determinant $\pm 1$
(the set of {\em unimodular} matrices).
Recall that any $n\times n$ matrix $[u_{ij}]$ with
$u_{ij}\!=\!0$ for all $i\!>\!j$ is called {\em upper triangular}.

Given any $M\!\in\!\Z^{n\times n}$, we then call an
identity of the form $UM=H$, with $H\!=\![h_{ij}]\!\in\!\Z^{n\times n}$
upper triangular and $U\!\in\!\gln(\Z)$, a {\em Hermite factorization}
of $M$. Also, if we have the following conditions in addition:
\begin{enumerate}
\item{$h_{ij}\!\geq\!0$ for all $i,j$.}
\item{\scalebox{.96}[1]{for all $i$, if $j$ is the smallest $j'$ such that
$h_{ij'}\!\neq\!0$ then $h_{ij}\!>\!h_{i'j}$ for all $i'\!\leq\!i$.}}
\end{enumerate}
then we call $H$ {\em \underline{the} Hermite normal form} of $M$. 

Also, given any identity of the form $UMV\!=\!S$ with $U,V\!\in\!\gln(\Z)$ 
and $S$ diagonal a {\em Smith factorization}. In particular, if 
$S\!=\![s_{i,j}]$ and we
require additionally that $s_{i,i}\!\geq\!0$ and $s_{i,i}|s_{i+1,i+1}$
for all $i\!\in\!\{1,\ldots,n\}$ (setting $s_{n+1,n+1}\!:=\!0$), then such a 
factorization for $M$ is unique and is called {\em \underline{the}} Smith 
factorization.

Finally, defining $x^A\:=\!(x^{a_{1,1}}_1
\cdots x^{a_{n,1}}_n,\ldots,x^{a_{1,n}}\cdots x^{a_{n,n}}_n)$, 
we call any map defined by\linebreak $x\mapsto x^A$ a {\em monomial change of
variables}. \dia 
\end{dfn} 
\begin{prop} 
\label{prop:mono} 
We have that $x^{AB}\!=\!(x^A)^B$ for any $A,B\!\in\!\Z^{n\times n}$. 
Also, for any field $K$, the map defined by $m(x)\!=\!x^U$,  
for any unimodular matrix $U\!\in\!\Z^{n\times n}$, is an automorphism of 
$(K^*)^n$. Finally, for any column vector $v\!\in\!\Zn$, 
the smallest valuation of an entry of $Uv$ is $k \Longleftrightarrow$ 
the smallest valuation of an entry of $v$ is $k$. \qed 
\end{prop} 
\begin{thm}
\label{thm:unimod}
\cite[Ch.\ 6 \& 8, pg.\ 128]{storjophd}
For any $A\!=\![a_{i,j}]\!\in\!\Z^{n\times n}$, the Hermite and Smith 
factorizations of $A$ can be computed within $O\!\left(n^{3.376} 
\log^2(n\max_{i,j}|a_{i,j}|)\right)$ bit operations.
Furthermore, the entries of all matrices in the Hermite and Smith 
factorizations have bit size $O(n\log(n\max_{i,j}|a_{i,j}|))$. \qed
\end{thm}
\begin{lemma}
\label{lemma:binomsys}
Following the notation of Definition \ref{dfn:smith} and  
Proposition \ref{prop:mono}, suppose $\det A\!\neq\!0$, 
$c_1,\ldots,c_n\!\in\!\Q^*_p$, $c\!:=\!(c_1,\ldots,c_n)$, 
$c'\!:=\!(c'_1,\ldots,c'_n)\!:=\!\left(\frac{c_1}{p^{\ord_p c_i}},\ldots,
\frac{c_n}{p^{\ord_p c_n}}\right)$, $L\!:=\!\max_i \ord_p s_{i,i}$, 
and let $v_1,\ldots,v_n$ be the columns of $V$. Then $x^A\!=\!c$ has a 
solution in $(\Q^*_p)^n$ iff (a) $(\ord_p c_1,\ldots,\ord_p c_n)v_i\!=\!0$ 
mod $s_{i,i}$ for all $i$ and (b) $x^A\!=\!c'$ has a solution 
in $((\Z/p^{2L+1})^*)^n$. 
In particular, the existence of a solution in $(\Q^*_p)^n$ for $x^A\!=\!c$  
can be decided in time polynomial in $n$ and $\log(n\max_{i,j}|a_{i,j}|)$. 
\end{lemma} 

\noindent 
{\bf Proof:} The necessity of Condition (a) follows immediately 
from Proposition \ref{prop:mono} upon observing that the valuations 
of the vector $x^A$ are exactly the entries of $[\ord_p x_1,
\ldots,\ord_p x_n]A$. Conversely, should Condition (a) hold, 
we can reduce to the case where $\ord_p c_i\!=\!0$ for all $i$. 
So let us assume the last condition. 

Observe now that $x^A\!=\!c$ iff $x^{AV}\!=\!c'$. Upon substituting 
$x\!:=\!y^U$, we see that the latter equation holds iff $y^{UAV}\!=\!c^V$. In 
other words, $y^S\!=\!c^V$. By Proposition \ref{prop:mono}, the last system 
has a solution in $(\Q^*_p)^n$ iff the first system does. By Corollary 
\ref{cor:binom} we thus see that Condition (b) is necessary and sufficient. 

To prove the asserted complexity bound, note that 
we can find $U$, $V$, and $S$ within the asserted time bound, 
thanks to Theorem \ref{thm:unimod}. Note also that by recursive 
squaring (and the observation that $\det A\!=\!\prod^n_{i=1} s_{i,i}$), 
we can find the $p$-parts of the $s_{i,i}$ and thus compute $L$ 
in polynomial-time. So then, applying Corollary \ref{cor:binom} $n$ times, we 
can decide in $\pp$ whether $y^S\!=\!c^V$ has a root in $(\Q^*_p)^n$.  \qed 

A final ingredient we will need is a method to turn roots  
of honest $n$-variate $(n+1)$-nomials on coordinate subspaces to roots  
in the algebraic torus. 
\begin{lemma} 
\label{lemma:toric} 
Suppose $c_0,\ldots,c_{k+1}\!\in\!\Q^*_p$, $a_1,\ldots,a_k\!\in\!\Z^k$ are 
linearly independent vectors, $\alpha\!:=\!(\alpha_1,\ldots,
\alpha_{k+1})\!\in\!\Z^{k+1}$ with $\alpha_{k+1}\!>\!0$, and 
$f(x)\!:=\!c_0+c_1x^{a_1}+\cdots+ c_kx^{a_k}+cx^\alpha$  
has a root in $(\Z_p\setminus\{0\})^k\times \{0\}$. Then $f$ 
has a non-degenerate root in $(\Z_p\setminus\{0\})^{k+1}$. \qed 
\end{lemma} 

\noindent
{\bf Proof:} Let $\zeta\!=\!(\zeta_1,\ldots,\zeta_k,0)\!\in\!(\Z_p\setminus 
\{0\})^k\times \{0\}$ be the stated root of $f$ and 
let $A$ denote the $k\times k$ matrix whose columns are $a_1,\ldots,a_k$. 
By Proposition \ref{prop:mat} we then have that $(\zeta_1,\ldots,\zeta_k)$ is 
a non-degenerate root of $\bar{f}(x)\!:=\!c_0+c_1x^{a_1}+\cdots+ c_kx^{a_k}$. 

To conclude, observe that $\frac{\partial \bar{f}}{\partial x_i}(\zeta_1,
\ldots,\zeta_k)\!=\!\frac{\partial f}{\partial x_i}(\zeta_1,\ldots,\zeta_k,
0)$ for all $i\!\in\!\{1,\ldots,k\}$. So $\zeta$ is a non-degenerate root of 
$f$. By the Implicit Function Theorem for analytic (i.e., $C^\infty$) 
functions over $\Q^n_p$ \cite[Thm.\ 7.4, pg.\ 237]{glockner}, there 
must then be a (non-degenerate) root $(\zeta'_1,\ldots,\zeta'_k,p^\ell)$ of 
$f$ for any sufficiently large $\ell\!\in\!\N$, with $\zeta'_i 
\longrightarrow \zeta_i$ for all $i\!\in\!\{1,\ldots,k\}$ as 
$\ell\!\longrightarrow +\infty$. Thus, we can find a root of $f$ in 
$(\Z_p\setminus\{0\})^{k+1}$. \qed 
\begin{rem} 
Note that Example \ref{ex:martin} from Section \ref{sub:rel} 
shows that the converse of Lemma \ref{lemma:toric} need not hold. 
On the other hand, over the real numbers, both the corresponding analogue of 
Lemma \ref{lemma:toric} and its converse hold \cite[Cor.\ 2.6]{brs}. \dia 
\end{rem} 

\medskip 
Henceforth, we will let $\bO$ 
denote the origin in whatever vector space we are working with. 
\begin{dfn} 
Suppose $K$ is a field, $c_0,\ldots,c_n\!\in\!K^*$, the 
vectors $a_0,\ldots,a_n\!\in\!\Zn$ are such that $a_1-a_0,\ldots,a_n-a_0$ 
are linearly independent, and 
$f(x)\!:=\!c_0x^{a_0}+c_1x^{a_1}+\cdots+c_nx^{a_n}$. We then 
call any sub-summand of the form $\bar{f}(x)\!=\!c_{i_1}x^{a_{i_1}}
+\cdots+c_{i_r}x^{a_{i_r}}$, with $\{i_1,\ldots,i_r\}$ of 
cardinality $r\!\geq\!1$, an {\em initial term polynomial} of $f$. \dia 
\end{dfn} 
\begin{rem} 
Note that setting any subset of variables equal to $0$ in 
$f$ --- with the resulting Laurent polynomial still well-defined and 
not identically $0$ --- 
results in an initial term polynomial of $f$. \dia 
\end{rem} 

\begin{cor} 
\label{cor:toric} 
Suppose $f\!\in\!\C_p\!\left[x^{\pm 1}_1,\ldots,x^{\pm 1}_n\right]$ 
has positive-dimensional Newton polytope with $\bO$  
as one of its vertices. Then $f$ has a root in $(\Q^*_p)^n 
\Longleftrightarrow$ some initial term polynomial of $f$ with 
at least $2$ terms has a root in $(\Q^*_p)^n$. 
\end{cor} 

\noindent 
{\bf Proof:} The ($\Longrightarrow$) direction is trivial since 
$f$ is an initial term polynomial by default. So let us 
focus on the ($\Longleftarrow$) direction. 

By assumption, we can then write $f(x)\!=\!c_0+c_1x^{a_1}+\cdots+c_nx^{a_n}$ 
with $c_0,\ldots,c_n\!\in\!\C^*_p$. Let $\zeta\!\in\!(\Q^*_p)^n$  
be a root of some initial term polynomial $\bar{f}$ of $f$. 
By Proposition \ref{prop:mono}, $\bar{f}(x)$ has a root in $(\Q^*_p)^n 
\Longleftrightarrow \bar{f}\!\left(x^U\right)$ has a root in $(\Q^*_p)^n$. So 
via the Hermite Factorization, we may assume that $f(x)\!=\!c_0
+c_1 x^{a_1}+\cdots+c_nx^{a_n}$ and the matrix $A$ whose
columns are $a_1,\ldots,a_n$ is upper-triangular. 
In other words, we may assume that $\bar{f}$ is independent of its 
last $n-r$ variables, for some $r\!\in\!\{1,\ldots,n-1\}$. 
So then, we may assume that $\zeta\!\in\!(\Q^*_p)^r\times \{0\}^{n-r}$ and 
$\bar{f}\!\in\!\C_p[x^{\pm 1}_1,\ldots,x^{\pm 1}_r]$. By multiplying certain 
rows of $A$ by $-1$ we can then clearly assume that 
$\zeta\!\in\!(\Z_p\setminus\{0\})^r
\times \{0\}^{n-r}$. By Lemma \ref{lemma:toric} (and 
induction) we then obtain that $f$ must have a root in 
$(\Q^*_p)^n$. \qed 

\subsection{The Proofs of Assertions (0) and (3) of Theorem \ref{thm:qp}} 
{\bf Assertion (0):} 
First note that the case $m\!\leq\!1$ is trivial:       
such a univariate $m$-nomial has no roots in $\Q_p$ iff  
it is a nonzero constant. 

The case $m\!=\!2$ then follows immediately  
from Corollary \ref{cor:binom}. \qed 

\medskip 
\noindent 
{\bf Assertion (3):}\\
{\bf Part (a):} 
First note that if $\zeta\!=\!(\zeta_1,\ldots,\zeta_n)\!\in\!\Q^n_p$ is a root 
of $f$ then all the exponents of $x_i$ in $f$ must be nonnegative for 
$\zeta_i\!=\!0$. We can then assume that, for all such $i$, some exponent of 
$x_i$ must be $0$. (Otherwise, $f$ would vanish on the entire hyperplane 
$\{y_i\!=\!0\}$, 
and the strict positivity of these exponents of $x_i$ in $f$ would be 
checkable a priori in quadratic time.) Note also that $\zeta$ being a root of 
$f$ is unaffected if we multiply $f$ by any power of $x_j$, provided 
$\zeta_j\!\neq\!0$. 

We can then clearly assume that $f$ has a nonzero constant term, 
write $f(x)\!=$\linebreak
$c_0+c_1x^{a_1}+\cdots+c_nx^{a_n}$ for some $c_0,\ldots,
c_n\!\in\!\Z\setminus\{0\}$, and let $A$ denote the matrix with 
columns $a_1,\ldots,a_n$. (Note also that enforcing our assumption that 
$f$ have a nonzero constant term induces at worst a factor of $2$ 
growth in absolute values of the entries of $A$.) By Corollary \ref{cor:toric} 
it then suffices to certify the existence of a root of $f$ in $(\Q^*_p)^n$. 

Set $L\!:=\!\max_i \ord_p(c_i)+\max_i \ord_p s_{i,i}+1$ 
where the $s_{i,i}$ denote the diagonal entries of the Smith Normal Form of 
$A$. Our certificate for $f$ having a root in $(\Q^*_p)^n$ will then be 
a root $\mu_0\!\in\!(\Z/p^{2L+1}\Z)^n\setminus\{\bO\}$ of 
the mod $p^{2L+1}$ reduction of 
$\bar{h}(x)\!:=\!\bar{g}(x^{\pm 1}_1,\ldots,x^{\pm 1}_n)$,  
for some choice of reciprocals, where 
$\bar{g}(x)\!:=\!x^{-a_i}\bar{f}(x)$ for some 
$i$, and $\bar{f}$ is an initial term polynomial of $f$ with at 
least $2$ terms. We will now show that $f$ has a root $\zeta\!\in\!(\Q^*_p)^n$ 
iff a certificate of the preceding form exists. 

To prove the ($\Longrightarrow$) direction, let us first 
clarify the choice of reciprocals in $\bar{g}(x^{\pm 1}_1,\ldots,x^{\pm 1}_n)$: 
we place an exponent of $-1$ for all $j$ where 
$\zeta_j\!\in\!\Q_p\setminus\Z_p$. Clearly then, with the preceding 
choice of reciprocals, $f(x^{\pm 1}_1,\ldots,x^{\pm 1})$ has a root 
$\mu\!\in\!(\Z_p\setminus\{0\})^n$. The choice of $i$ to 
define $\bar{h}(x)$ is also simple to pin down:  
pick any $i$ with $\ord_p(\mu^{a_i})$ minimal. The roots of 
$h(x)\!:=\!x^{-a_i}f(x^{\pm 1}_1,\ldots,x^{\pm 1}_n)$ in 
$(\Q^*_p)^n$ are clearly independent of $i$. 

To clarify the choice of $\bar{f}$ let us first write 
$h(x)\!:=\!\gamma_0+\gamma_1x^{\alpha_1}+\cdots+\gamma_nx^{\alpha_n}$.  
The $\gamma_i$ are then a re-ordering of the $c_i$, the 
$\alpha_i$ are differences of columns of $A$, and the matrix $A'$ 
with columns $\alpha_1,\ldots,\alpha_n$ is non-singular and has entries 
no larger in absolute value than twice those of $A$.  
We also have that $\ord_p(\mu^{\alpha_i})\!\geq\!0$ for all $i$ by 
construction. Moreover, by the ultrametric property (applied to the
sum $\gamma_0+(c_1\mu^{\alpha_1}+\cdots+\gamma_n\mu^{\alpha_n})$), the 
root $\mu$ of $h$ must satisfy $\ord_p(\gamma_i\mu^{a_i})\!\leq\!\ord_p 
\gamma_0\!\leq\!\max_k \ord_p c_k\!\leq\!L$ for some $i$. (Otherwise 
$\ord_p h(\mu)=\ord_p \gamma_0\!<\!+\infty$).
By Propositions \ref{prop:mat} and \ref{prop:mono}, and 
the Smith factorization of the matrix $A'$, 
we must then have $\ord_p h_j (\mu)\!\leq\!\ord_p(\gamma_0)+\max_i 
\ord_p(2s_{i,i})\!\leq\!L\!=\!O(\size(f))$ for some $j$. 

Clearly then, there are 
$u_{i_1},\ldots,u_{i_r}\!\in\!\Z_p\setminus\{0\}$ with $r\!\geq\!1$, 
$L\!\geq\!\ord_p u_{i_j}\!\geq\!\ord_p\gamma_{i_j}$ for all $j$, 
$\gamma_0+u_{i_1}+\cdots+u_{i_r}\!=\!0$,  
and $(\mu^{\alpha_{i_1}},\ldots,\mu^{\alpha_{i_r}})\!=\!
\left(\frac{u_{i_1}}{c_{i_1}},\ldots,\frac{u_{i_r}}{c_{i_r}}\right)$. So 
define $\bar{f}$ to be the sum of terms of $f$ corresponding to picking 
the $i_1,\ldots,i_r$ terms of $h$. By Lemma \ref{lemma:binomsys}, $\mu$ then 
has a well-defined mod $p^{2L+1}$ reduction 
$\mu_0\!\in\!(\Z/p^{2L+1}\Z)^n\setminus\{\bO\}$ 
that is a root of the mod $p^{2L+1}$ reduction of $\bar{h}$. 
So the ($\Longrightarrow$) direction is proved.  

To prove the ($\Longleftarrow$) direction, let us suppose that 
the mod $p^{2L+1}$ reduction of 
$\bar{h}(x)\!:=\!\bar{g}(x^{\pm 1}_1,\ldots,x^{\pm 1})$ 
has a root $\mu_0\!\in\!(\Z/p^{2L+1}\Z)^n\setminus\{\bO\}$ for some 
choice of signs, some choice of $i$, and some choice of 
initial term polynomial $\bar{f}$ of $f$ so that 
$\bar{g}(x)\!=\!x^{-a_i}\bar{f}(x)$. 
Writing $\bar{h}(x)\!=\!\gamma_0+\gamma_{i_1}x^{\alpha_{i_1}}+\cdots+
\gamma_{i_r}x^{\alpha_{i_r}}$ as before, it is clear that 
$\ord_p(\gamma_i \mu^{\alpha_i})\!\leq\!
\ord_p \gamma_0$ for some $i$ by the ultrametric inequality. 
So then, by Proposition \ref{prop:mat}, $\ord_p \bar{h}'(\mu)\!\leq\!L$, and 
then by Hensel's Lemma, $\bar{h}$ has a root 
$\mu'\!\in\!\Zn_p\setminus\{\bO\}$. By Corollary \ref{cor:toric}, 
$h(x)\!:=\!\gamma_0+\gamma_1x^{\alpha_1}+\cdots+\gamma_nx^{\alpha_n}$ must 
then have a root $\mu\!\in\!(\Z_p\setminus\{\bO\})^n$. So by the definition 
of $h$, it is then clear that defining $\zeta_i\!=\!\mu^{\pm 1}_i$ for a 
suitable choice of signs, $\zeta\!:=\!(\zeta_1,\ldots,\zeta_n)$ is a 
root of $f$. \qed 

\noindent 
{\bf Part (b):} Since the Legendre symbol $\left(\frac{a}{p}
\right)$ can be evaluated within $O((\log a)(\log p))$ bit operations 
\cite[Thm.\ 5.9.3, pg.\ 113]{bs}, the criteria from Theorem \ref{thm:quad} 
can clearly be checked in time polynomial in $\size(f)$. So we are done. \qed 

\noindent 
{\bf Part (c):} 
Via the Smith Normal Form, Proposition \ref{prop:mono}, and 
Corollary \ref{cor:toric}, we can reduce to the  
special case detailed in Theorem \ref{thm:weil}, 
i.e., we may assume that we have an instance of the form 
$f(x)\!=\!c_0+c_1x^{d_1}_1+\cdots+c_nx^{d_n}$ with $d_1,\ldots,d_n\!\in\!\N$,  
and thus $n!V_f\!=\!\prod^n_{i=1} d_i\!>\!\prod^n_{i=1}(\gcd(d_i,p-1)-1)$. 

By the succinct certificates we used to prove Part (a), we see that 
the existence of a root of $f$ in $\Q^n_p$ is implied by the existence of a 
root of $f$ in $\F^n_p$ if $\ord_p|c_0|\!=\ldots=\!\ord_p|c_n|\!=\!
\ord_p(n!V_f)\!=\!0$. By Theorem \ref{thm:weil}, a root for $f$ in $\F^n_p$ 
is guaranteed if $n\!\geq\!2$, $p$ does not divide any $c_i$, and 
$p\!\geq\!(n!V_f)^{2/(n-1)}$. So we are done. \qed 

\section{Discriminants, $p$-adic Newton Polygons, and\\ Assertion (2)} 
\label{sec:disc} 
The intuition behind the speed-up of Assertion (2) is that the hardness of 
instances of $\fqp(\Z[x_1]\times \Pro)$ is governed by numerical conditioning, 
quite similar to the sense long known in
numerical linear algebra (and extended more recently to
real feasibility \cite{cuckersmale}). More concretely, the classical fact that
Newton iteration converges more quickly for a root $\zeta\!\in\!\C$ of $f$
with $f'(\zeta)$ having large norm (i.e., a {\em well-conditioned} root)
persists over $\Q_p$.

To prepare for our next proof, let us first clarify the statement about natural 
density $0$ in Assertion (2) of Theorem \ref{thm:qp}. 
\begin{dfn}
\label{dfn:basic}
Letting $\#$ denote set cardinality, we say that
$S\!\subseteq\!\Pro$ {\em has (natural) density $\mu$} iff
$\lim\limits_{t\rightarrow \infty}\frac{\#S\cap\{1,\ldots,t\}}
{\#\Pro\cap\{1,\ldots,t\}}\!=\!\mu$. \dia
\end{dfn}

\noindent 
Now let $(\Z\times (\N\cup \{0\}))^\infty$ denote the set of all infinite
sequences of pairs $((c_i,a_i))^\infty_{i=1}$
with $c_i\!=\!a_i\!=\!0$ for $i$ sufficiently large.
Note then that $\Z[x_1]$ admits a natural embedding
into $(\Z\times (\N\cup \{0\}))^\infty$ by
considering coefficient-exponent pairs in order of increasing
exponents, e.g., $a+bx^{99}+x^{2001} \mapsto
((a,0),(b,99),(1,2001),(0,0),(0,0),\ldots)$. Then natural density
for a set of pairs $\cI\!\subseteq\!\Z[x_1]\times \Pro$ then simply means
the corresponding natural density within $(\Z\times (\N\cup \{0\}))^\infty
\times \Pro$. 

The exceptional set to Assertion (2) can be made more precise once 
one introduces the $\cA$-discriminant. But first we must introduce the 
resultant and some quantitative estimates. 
\begin{dfn}
\label{dfn:syl}
(See, e.g., \cite[Ch.\ 12, Sec.\ 1, pp.\ 397--402]{gkz94}.)
Suppose\linebreak $f(x_1)\!=\!a_0+\cdots+a_dx^{d}_1$ and
$g(x_1)\!=\!b_0+\cdots+b_{d'}x^{d'}_1$ are polynomials with
indeterminate coefficients. We define their {\em Sylvester matrix}
to be the $(d+d')\times (d+d')$ matrix

\noindent
\mbox{}\hfill \scalebox{1}[1]{$\cS_{(d,d')}(f,g)\!:=\!\begin{bmatrix}
a_0 & \cdots & a_d    & 0       & \cdots & 0 \\
   & \ddots &  & & \ddots &  \\
0   & \cdots & 0 & a_0 & \cdots & a_d \\
b_0 & \cdots & b_{d'}    & 0       & \cdots & 0 \\
  & \ddots &  &  & \ddots &  \\
0   & \cdots & 0 & b_0 & \cdots & b_{d'} 
\end{bmatrix}
\begin{matrix}
\\
\left. \rule{0cm}{.9cm}\right\}
d' \text{ rows}\\
\left. \rule{0cm}{.9cm}\right\}
d \text{ rows} \\
\\
\end{matrix}$}\hfill\mbox{}\\
and their {\em Sylvester resultant} to be
$\cR_{(d,d')}(f,g)\!:=\!\det \cS_{(d,d')}(f,g)$. \dia 
\end{dfn}
\begin{lemma}
\label{lemma:syl}
Following the notation of Definition \ref{dfn:syl},
assume $f,g\!\in\!K[x_1]$ for some field $K$, and that
$a_d$ and $b_{d'}$ are not both $0$. Then $f\!=\!g\!=\!0$ has a root in the
algebraic closure of $K$ iff $\cR_{(d,d')}(f,g)\!=\!0$. More
generally, we
have $\cR_{(d,d')}(f,g)\!=\!a^{d'}_d \!\! \prod\limits_{f(\zeta)=0} g(\zeta)$
where the product counts multiplicity. Finally, if we assume further that $f$
and $g$ have complex coefficients of absolute value $\leq\!H$, and $f$
(resp.\ $g$) has exactly $m$ (resp.\ $m'$) monomial terms, then
$|\cR_{(d,d')}(f,g)|\!\leq\!m^{d'/2}m'^{d/2}H^{d+d'}$. \qed
\end{lemma}
\noindent
The first $2$ assertions are classical (see, e.g.,
\cite[Ch.\ 12, Sec.\ 1, pp.\ 397--402]{gkz94} and \cite[pg.\ 9]{rs}).
The last assertion follows easily from
Hadamard's Inequality (see, e.g., \cite[Thm.\ 1, pg.\ 259]{mignotte}).

We are now ready to introduce discriminants. 
\begin{dfn}
\label{dfn:adisc} For any field $K$, write any $f\!\in\!K[x_1]$ as
$f(x_1)\!=\!\sum^m_{i=1}c_ix^{a_i}_1$
with $0\!\leq\!a_1\!<\cdots<\!a_m$. Letting $\cA\!=\!\{a_1,\ldots,a_m\}$,
we then define
the {\em $\cA$-discriminant} of $f$, $\Delta_\cA(f)$, to be\\
\mbox{}\hfill $\cR_{(\bar{a}_m,\bar{a}_m-\bar{a}_2)}\left.
\left(\bar{f},\left.\frac{\partial\bar{f}}{\partial x_1}
\right/x^{\bar{a}_2-1}_1\right)\right/
c^{\bar{a}_m-\bar{a}_{m-1}}_m$,
\hfill\mbox{}\\
where $\bar{a}_i\!:=\!(a_i-a_1)/g$ for all $i$, $\bar{f}(x_1)\!:=\!
\sum^m_{i=1}c_ix^{\bar{a}_i}_1$, and $g\!:=\!\gcd(a_2-a_1,\ldots,
a_m-a_1)$ (see also \cite[Ch.\ 12, pp.\ 403--408]{gkz94}). Finally, if
$c_i\!\neq\!0$ for all $i$, then we call $\supp(f)\!:=\!\{a_1,\ldots,a_m\}$
the {\em support} of $f$. \dia
\end{dfn}
\begin{rem}
Note that when $\cA\!=\!\{0,\ldots,d\}$ we have
$\Delta_\cA(f)\!=\!\cR_{(d,d-1)}(f,f')/c_d$, i.e.,
for dense polynomials, the $\cA$-discriminant agrees with the
classical discriminant \dia
\end{rem}

The claim of natural density $0$ in Assertion (2) of Theorem \ref{thm:qp} 
can then be made explicit as follows. 
\begin{cor}
\label{cor:lots}
For any subset $\cA\!=\!\{a_1,\ldots,a_m\}\!\subset\!\N\cup\{0\}$
with $0\!=\!a_1\!<\cdots<\!a_m$, let $T_\cA$ denote the family of pairs
$(f,p)\!\in\!\Z[x_1]\times \Pro$ with $f(x_1)\!=\!\sum^m_{i=1}c_ix^{a_i}_1$
and let $T^*_\cA$ denote the subset of $T_\cA$
consisting of those pairs $(f,p)$ with
$p\not\!|\Delta_\cA(f)$. Also let $T_\cA(H)$ (resp.\ $T^*_\cA(H)$)
denote those pairs $(f,p)$ in $T_\cA$ (resp.\ $T^*_\cA$) where
$|c_i|\!\leq\!H$ for all $i\!\in\![m]$ and $p\!\leq\!H$. Finally,
let $d\!:=\!a_m/\gcd(a_2,\ldots,a_m)$. Then for all
$H\!\geq\!17$ we have\\
\mbox{}\hfill
$\frac{\#T^*_\cA(H)}{\#T_\cA(H)}\!\geq\!\left(1-\frac{(2d-1)m}{2H+1}\right)
\left(1-\frac{1+(2d-1)\log(mH)\log H}{H}\right)$.\hfill \mbox{}
\end{cor}

\noindent
In particular, we will see in the proof of Assertion (2) of
Theorem \ref{thm:qp} that the exceptional set $\cE$ is merely the
complement of the union $\bigcup_{\cA} \cT^*_\cA$ as $\cA$ ranges over
all finite subsets of $\N\cup\{0\}$. Our corollary above is proved in Section
\ref{sub:lots}.

Another bit of background we'll need to prove Assertion (2) 
of Theorem \ref{thm:qp} is some arithmetic tropicalia. 
\begin{dfn} Given any polynomial
$f(x_1)\!:=\!\sum^m_{i=1}c_ix^{a_i}_1$
$\in\!\Z[x_1]$, we define its {\em $p$-adic Newton polygon},
$\newt_p(f)$, to be the convex hull of
the points $\{(a_i,\ord_p c_i)\; | \; i\!\in\!\{1,\ldots,m\}\}$.
Also, a face of a polygon $P\!\subset\!\R^2$ is called {\em lower} iff it has
an inner normal with positive last coordinate, and the {\em lower hull} of
$P$ is simply the union of all its lower edges. Finally,
the polynomial associated to summing the terms of $f$ corresponding to
points of the form $(a_i,\ord_p c_i)$ lying on a lower
face of $\newt_p(f)$
is called a {\em ($p$-adic) lower polynomial}. \dia
\end{dfn}
\begin{ex}
\label{ex:newt}
For $f(x_1):=36 -8868x_1 +29305x^2_1 -35310x^3_1 +18240x^4_1
-3646x^5_1+243x^6_1$,\linebreak

\vspace{-.8cm}
\noindent
\begin{minipage}[t]{3in}
\vspace{0pt}
the polygon $\newt_3(f)$ has exactly $3$ lower edges and
can easily be verified to resemble the
illustration to the right. The polynomial $f$ thus has
exactly $2$ lower binomials, and $1$ lower trinomial. \dia
\end{minipage}
\hfill\raisebox{-2.5cm}{\epsfig{file=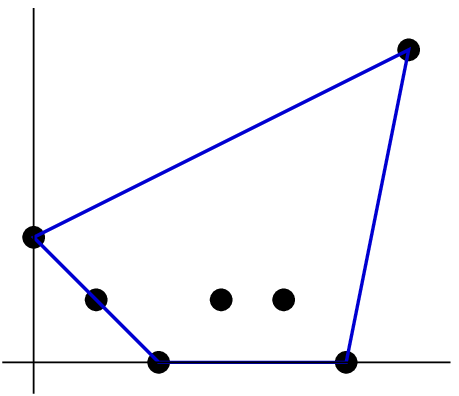,height=.95in}}
\hfill\mbox{}
\end{ex}

A remarkable fact true over $\C_p$ but false over $\C$ is that 
the norms of roots can be determined completely combinatorially. 
\begin{lemma}
\label{lemma:newt}
(See, e.g., \cite[Ch.\ 6, sec.\ 1.6]{robert}.)
The number of roots of $f$ in $\C_p$ with valuation $v$, counting
multiplicities, is {\em exactly} the horizontal length of the lower face of
$\newt_p(f)$ with inner normal $(v,1)$. \qed
\end{lemma}
\begin{ex}
In Example \ref{ex:newt} earlier, note that the $3$ lower
edges have respective horizontal lengths $2$, $3$, and $1$,
and inner normals $(1,1)$, $(0,1)$, and $(-5,1)$. Lemma \ref{lemma:newt}
then tells us that $f$ has exactly $6$ roots in $\C_3$:
$2$ with $3$-adic valuation $1$, $3$ with $3$-adic valuation $0$, and $1$
with $3$-adic valuation $-5$. Indeed, one can check that the roots
of $f$ are exactly $6$, $1$, and $\frac{1}{243}$, with respective
multiplicities $2$, $3$, and $1$. \dia
\end{ex}

\subsection{The Proof of Assertion (2) of Theorem \ref{thm:qp}}  
\label{sub:proof2} 
The existence of $0$ as a root is clearly checkable in constant time
so we may again assume that $f$ is not divisible by $x_1$.
Via the reciprocal polynomial $f^*(x_1)\!:=\!x_1^{\deg f}f(1/x_1)$, it is 
then enough to 
show that, for most $f$, having a root in $\Z_p$ admits a succinct certificate.
As observed in the proof of Assertion (2), $\newt_p(f)$ can
be computed in polynomial-time. Since
$\ord_pc_i\!\leq\!\log_pc_i\!\leq\!\size(c_i)$, note also that
that every root $\zeta\!\in\!\C_p$ of $f$
satisfies $|\ord_p\zeta|\!\leq\!2\max_i\size(c_i)\!\leq\!2\size(f)\!<\!
2\size_p(f)$.

Since $\ord_p(\Z_p)\!=\!\N\cup\{0\}$, we can clearly assume that
$\newt_p(f)$ has an edge with non-positive integral slope, for otherwise $f$
would have no roots in $\Z_p$. Letting $g(x_1)\!:=\!f'(x_1)/x_1^{a_1-1}$,
and $\zeta\!\in\!\Z_p$ be any $p$-adic integer root of $f$, note then that
$\ord_p f'(\zeta)\!=\!(a_1-1)\ord_p(\zeta)+\ord_p g(\zeta)$. Note also that
$\Delta_\cA(f)\!=\!\res_{a_m,a_m-a_1}(f,g)$ so if $p\not|\Delta_\cA(f)$
then $f$ and $g$ have no common roots in the algebraic closure of
$\F_p$, by Lemma \ref{lemma:syl}. In particular,
$p\!\!\not|\Delta_\cA(f)\Longrightarrow g(\zeta)\!\not
\equiv\!0 \; \mod \; p$; and thus $p\!\!\not\!|\Delta_\cA(f,g)\Longrightarrow
\ord_p f'(\zeta)\!=\!(a_1-1)\ord_p(\zeta)$. Furthermore, by the convexity
of the lower hull of $\newt_p(f)$, it is clear that $\ord_p(\zeta)\!\leq\!
\frac{\ord_p c_0 -\ord_p c_i}{a_i}$ where $(a_i,\ord_p c_i)$ is
the rightmost vertex of the lower edge of $\newt_p(f)$ with
least (non-positive and integral) slope. Clearly then,
$\ord_p(\zeta)\!\leq\!\frac{2\max_i \log_p|c_i|}{a_1}$.
So $p\not\!|\Delta_\cA(f)\Longrightarrow \ord_p f'(\zeta)\!\leq\!2\size(f)$.

Our fraction of inputs admitting a succinct certificate will then
correspond precisely to those $(f,p)$ such that $p\!\!\not\!|\Delta_\cA(f)$.
In particular, let us define $\cE$ to be the union of all pairs
$(f,p)$ such that $p|\Delta_\cA(f)$, as $\cA$ ranges over all finite
subsets of $\N\cup\{0\}$. It is then easily checked that $\cE$ is a
countable union of hypersurfaces, and the density $0$ statement 
follows immediately from Corollary \ref{cor:lots}. 

Now fix $\ell\!=\!4\size(f)+1$. Clearly then,
by Hensel's Lemma, for any $(f,p)\!\in\!(\Z[x_1]\times \Pro)\setminus \cE$,
$f$ has a root $\zeta\!\in\!\Z_p \Longleftrightarrow f$ has a root
$\zeta_0\!\in\!\Z/p^\ell\Z$. Since
$\log(p^\ell)\!=\!O(\size(f)\log p)\!=\!O(\size_p(f)^2)$,
and since arithmetic in $\Z/p^\ell\Z$
can be done in time polynomial in $\log(p^\ell)$ \cite[Ch.\ 5]{bs}, we have
thus at last found our desired certificate: a root
$\zeta_0\!\in\!(\Z/p^\ell\Z)^*$ of $f$ with $\ell\!=\!4\size(f)+1$. \qed

\section{Degenerate Trinomials, Linear Forms in p-adic\\ 
Logarithms, and Assertion (1)}  
\label{sec:lin} 

We will first need to recall the concept of a {\em gcd-free basis}. In
essence, a gcd-free basis is nearly as powerful as factorization into primes,
but is far easier to compute.
\begin{dfn}
\label{dfn:gcd}
\cite[Sec.\ 8.4]{bs} For any subset $\{\alpha_1,\ldots,\alpha_N\}$
$\subset\!\N$, a {\em gcd-free basis for $\{\alpha_1,\ldots,\alpha_N\}$} is a
pair of sets $\left(\{\gamma_i\}^\eta_{i=1}, \{e_{ij}\}_{(i,j)\in [N]\times 
[\eta]}\right)$ such that
(1) $\gcd(\gamma_i,\gamma_j)\!=\!1$ for all $i\!\neq\!j$, and
(2) $\alpha_i\!=\!\prod^{\eta}_{j=1} \gamma^{e_{ij}}_j$ for all $i$. \dia
\end{dfn}
\begin{thm}
\label{thm:gcd}
Following the notation of Definition \ref{dfn:gcd}, we
can compute a gcd-free basis for $\{\alpha_1,\ldots,\alpha_N\}$
(with $\eta$ linear in $N+\max_i\log \alpha_i$)
in time linear in $N+\max_i\log^2\alpha_i$.
In particular, if $u_1,\ldots,u_N\!\in\!\Z$ then we can decide
$\alpha^{u_1}_1\cdots \alpha^{u_N}_N\!\stackrel{?}{=}\!1$
in time linear in\linebreak
$N+\left(\max_i\log(\alpha_i)+\max_i\log(u_i)\right)^2$. \qed
\end{thm}

\noindent
The first assertion of Theorem \ref{thm:gcd} follows immediately from
\cite[Thm.\ 4.8.7, Sec.\ 4.8] {bs} and the naive bounds for the complexity of
integer multiplication. The second assertion then follows immediately
by checking whether the linear combinations $\sum^N_{i=1} e_{ij}u_i$
are all $0$ or not.

We now make some final observations about the roots of trinomials 
before proving\linebreak Assertion (1) of Theorem \ref{thm:qp}. 
\begin{cor}
\label{cor:tri} 
Suppose $f(x_1)\!=\!c_1+c_2x^{a_2}_1+c_3x^{a_3}_1\!\in\!\cF_{1,3}$,  
$\cA\!:=\!\{0,a_2,a_3\}$, $0\!<\!a_2\!<\!a_3$, $a_3\!\geq\!3$, and 
$\gcd(a_2,a_3)\!=\!1$. Then: 

\smallskip 
\noindent
(0) $\Delta_\cA(f)=(a_3-a_2)^{a_3-a_2}a^{a_2}_2c^{a_3}_2-
(-a_3)^{a_3}c^{a_3-a_2}_1c^{a_2}_3$.\\ 
(1) $\Delta_\cA(f)\!\neq\!0\Longleftrightarrow f$ has no degenerate
roots. In which case, we also have\\
\mbox{}\hfill$\Delta_\cA(f)\!=\frac{(-1)^{a_3}c^{a_2-1}_3}{c^{a_2-1}_1}
\!\!\!\!\prod\limits_{f(\zeta)=0}f'(\zeta)$.\hfill\mbox{}\\
(2) Deciding whether $f$ has a degenerate root in $\C_p$ can be  
done in time polynomial in\linebreak 
\mbox{}\hspace{.75cm}$\size_p(f)$.\\  
(3) If $f$ has a degenerate root $\zeta\!\in\!\C^*_p$ then  
$(\zeta^{a_2},\zeta^{a_3})\!=\!\frac{c_1}{a_3-a_2} 
\left(-\frac{a_3}{c_2},\frac{a_2}{c_3}\right)$. 
In particular, such a\linebreak 
\mbox{}\hspace{.7cm}$\zeta$ is unique and lies in $\Q$.\\ 
(4) The polynomial $q(x_1)\!:=\!(a_3-a_2)
-a_3x^{a_2}_1+a_2x^{a_3}_1$ has $1$
as its unique degenerate root and\linebreak 
\mbox{}\hspace{.7cm}satisfies 
$\Delta_{\{0,\ldots,a_3-2\}}\!\left(\frac{q(x_1)}{(x_1-1)^2}\right)=
\pm (a_2a_3(a_3-a_2))^{a_3+O(1)}$. 
\end{cor} 

\medskip 
\noindent 
{\bf Proof of Corollary \ref{cor:tri}:} 

\smallskip
\noindent 
{\bf Part (0):} \cite[Prop.\ 1.8, pg.\ 274]{gkz94}. \qed 

\smallskip
\noindent 
{\bf Part (1):} The first assertion follows directly from 
Definition \ref{dfn:adisc} and the vanishing criterion for 
$\res_{(a_3,a_3-a_2)}$ from 
Lemma \ref{lemma:syl}. To prove the second assertion, observe that 
the product formula from Lemma \ref{lemma:syl} implies that\\ 
$\Delta_\cA(f)  =  c^{a_3-a_2}_3 \!\! \left.\left(\prod_{f(\zeta)=0} 
\frac{f'(\zeta)}{\zeta^{a_2-1}}\right)\right/c^{a_3-a_2}_3=
(-1)^{a_3}\left.\left(\prod_{f(\zeta)=0}
f'(\zeta)\right)\right/(c_1/c_3)^{a_2-1}$. \qed  

\smallskip
\noindent
{\bf Part (2):} From Part (1) it suffices to detect the vanishing
of $\Delta_\cA(f)$. However, while Part (0) implies that
one can evaluate $\Delta_\cA(f)$ with a small number of arithmetic operations,
the bit-size of $\Delta_\cA(f)$ can be quite large. Nevertheless, we can
decide within time polynomial in $\size(f)$ whether
these particular $\Delta_\cA(f)$ vanish for integer $c_i$ via gcd-free
bases (invoking Theorem \ref{thm:gcd}). \qed

\smallskip
\noindent 
{\bf Part (3):} 
It is easily checked that if $\zeta\!\in\!\C_p$ is a 
degenerate root of $f$ 
then the vector $\left[c_1,c_2\zeta^{a_2},c_3\zeta^{a_3}\right]$ must be a 
right null vector for the matrix $M\!:=\!\begin{bmatrix}1 & 1 & 1 \\ 
0 & a_2 & a_3
\end{bmatrix}$. Since \linebreak 
$[a_3-a_2,-a_3,a_2]$ is clearly a right null 
vector for $M$, $\left[c_1,c_2\zeta^{a_2},c_3\zeta^{a_3}\right]$ must then be 
a mutiple of $[a_3-a_2,-a_3,a_2]$. Via the extended Euclidean algorithm
\cite[Sec.\ 4.3]{bs}, we can
then find $A$ and $B$ (also of size
polynomial in $\size(f)$) with $Aa_2+Ba_3\!=\!1$. So then we 
obtain that\\
\scalebox{.92}[1]{$\left(\frac{c_2\zeta^{a_2}}{c_1}\right)^A
\left(\frac{c_3\zeta^{a_3}}
{c_1}\right)^B\!=\!\frac{c^A_2c^B_3}{c^{A+B}_1}\zeta\!=\!
\left(\frac{-a_3}{a_3-a_2}\right)^A\left(\frac{a_2}
{a_3-a_2}\right)^B$.} \qed 

\smallskip
\noindent 
{\bf Part (4):} That $1$ is a root of $q$ is obvious. 
Uniqueness follows directly from Part (3) and our assumption 
that $\gcd(a_2,a_3)\!=\!1$. To prove the final assertion, first note 
that a routine long division reveals that $\frac{q(x)}{(x-1)^2}$ 
has coefficients rising by one arithmetic progression and then 
falling by another. Explicitly,\\ 
\mbox{}\hfill 
\scalebox{1}[1]{$\displaystyle{\frac{q(x)}{(x-1)^2}=\left(\sum^{a_2-1}_{i=1} 
(a_3-a_2)ix^{i-1} \right)+\left(\sum^{a_3-a_2}_{i=1} (a_3-a_2+1-i)
a_2x^{a_2-2+i}\right)}$.}\hfill \mbox{}\\ 
  
Definition \ref{dfn:syl} then implies that $\Delta_{\{0,
\ldots,a_3-2\}}\!\left( \frac{q(x_1)}{(x-1)^2}\right)$ is exactly 
$\frac{1}{a_2}$ times the determinant of the following 
quasi-Toeplitz matrix which we will call $\cM$: \\
\scalebox{.69}[.9]{$\begin{bmatrix}
a_3-a_2 & 2(a_3-a_2)   & \cdots & (a_2-1)(a_3-a_2) & (a_3-a_2)a_2 & \cdots  
 & 2a_2 & a_2 & 0      & \cdots & 0 \\
    & \ddots & \ddots &              & \ddots            & \ddots 
 &                  & \ddots  & \ddots &        &    \\ 
1\cdot 2\cdot (a_3-a_2) & 2\cdot3\cdot (a_3-a_2) & \cdots & 
 (a_2-2)(a_2-1)(a_3-a_2) & (a_2-1)(a_3-a_2)a_2 & \cdots 
 & (a_3-2)\cdot 1\cdot a_2 & 0 & \cdots 
 & & 0\\ 
    & \ddots & \ddots &              & \ddots            & \ddots 
 &                  & \ddots  & \ddots &        &    \\ 
\end{bmatrix}$},\\
where there are exactly $a_3-3$ (resp. $a_3-2$) shifts of the first 
(resp.\ second) detailed row. Letting $f(x)\!:=\!\frac{q(x)}{(x-1)^2}$, note 
in particular that the entries of 
the first $a_3-3$ (resp.\ last $a_3-2$) rows correspond to the coefficients of 
$x^if(x)$ (resp.\ $x^if'(x)$) for $i\!\in\!\{0,\ldots,a_3-4\}$ 
(resp.\ $i\!\in\!\{0,\ldots,a_3-3\}$). 
We can clearly replace any polynomial by itself plus a linear combination of 
the others and rebuild our matrix $\cM$ with these new polynomials, leaving 
$\det \cM$ unchanged (thanks to invariance under elementary row operations). 
So let us now look for useful linear combinations of $x^if$ and $x^jf'$. 

Observe that\\
\mbox{}\hfill
$\displaystyle{\frac{q(x)}{x-1}=\sum_{i=0}^{a_2-1}(a_2-a_3)x^i+
\sum_{i=a_2}^{a_3-1}a_2x^i}$ \hfill and \hfill 
$\displaystyle{\frac{q'(x)}{x-1}=a_2a_3x^{a_2-1}+\cdots+a_2a_3x^{a_3-2}},$
\hfill\mbox{}\\
so\\
\mbox{}\hfill
$\displaystyle{\frac{q(x)}{(x-1)}-\frac{1}{a_3}\frac{xq'(x)}{(x-1)}
=\sum_{i=0}^{a_2-1}(a_2-a_3)x^i}.$
\hfill\mbox{}\\ 
Since $(x-1)f(x)\!=\!\frac{q(x)}{x-1}$ it would thus be useful to 
obtain $\frac{q'(x)}{x-1}$ as a polynomial linear combination of $f$ and 
$f'$. Toward this end, observe that
\begin{align*}
  xf'-&f'+2f=(x-1)f'+2f\\
  &=\frac{(x-1)^2q'-2(x-1)q}{(x-1)^3}
  +\frac{2(x-1)q}{(x-1)^3}\\
  &=\frac{(x-1)^2q'}{(x-1)^3}=\frac{q'}{(x-1)}.
\end{align*} 
It is then prudent to replace each $x^i f$ row with the 
coefficients of\\
\mbox{}\hfill $x^i\left(f+\left(\frac{2}{a_3}-1\right)xf-\frac{x}{a_3}f'
+\frac{x^2}{a_3}f'\right)$, \hfill\mbox{}\\
for $i\!\in\!\{0,\ldots,a_3-5\}$. 
There are $a_3-4$ such new rows, each divisible by $a_3-a_2$, so 
$(a_3-a_2)^{a_3-4}$ divides $\det \cM$. Similarly, we can replace each 
$x^i f'$ row with the coefficients of $x^i(f'-xf'-2f)$, for 
$i\!\in\!\{0,\ldots,a_3-4\}$. Each of these polynomials is divisible by 
$a_2a_3$. There are $a_3-3$ of these rows --- and they are distinct from the 
other $a_3-4$ rows we modified earlier --- so $(a_2a_3)^{a_3-3}$ also divides 
$\det \cM$.

We are thus left with showing that the matrix whose rows 
correspond to the coefficient vectors of the polynomials\\ 
\mbox{}\hfill $\frac{x^{a_2}-1}{x-1}, \ldots,   
x^{a_3-5}\frac{x^{a_2}-1}{x-1},x^{a_3-4}f,x^{a_2-1}\frac{x^{a_3-a_2}-1}
{x-1},\ldots,x^{a_2+a_3-5}\frac{x^{a_3-a_2}-1}{x-1},x^{a_3-3}f'$, \hfill 
\mbox{}\\ 
has determinant $\pm (a_2a_3)^{O(1)}$. Roughly, our last matrix has 
the following form:\\
\mbox{}\hfill 
$\begin{bmatrix}
1 & \cdots  & 1      &       &        &  &  & \\
  & 1       & \cdots & 1     &        &  &  & \\  
  &         & \ddots &       & \ddots &  &  &  \\
  &         &        &  1    & \cdots & 1 &  & \\  
 &   &        &    & a_3-a_2 & \cdots & \cdots & a_2\\     
 &  & 1      & \cdots          & 1 &       &        & \\
  &         &  &   &  \ddots  & \ddots & &  \\
 &         &        &      &      & 1      & \cdots & 1\\  
 &  &        &     & 2(a_3-a_2) & \cdots & \cdots & (a_3-2)a_2\\     
\end{bmatrix}$\hfill \mbox{}\\  
Via a simple sequence of $O(a_3)$ elementary row and column operations, 
restricted 
to subtractions of a column from another column and subtractions of a 
row from another row, 
we can then reduce our matrix to a $(2a_3-5)\times (2a_3-5)$ 
permutation matrix with the ${a_3}^{\text{\underline{rd}}}$ row and 
$(2a_3-5)\thth$ row resembling  
the corresponding rows above. In particular, these $2$ new rows have 
entries at worst $O(a_3)$ times larger than before. Clearly then, 
our final determinant is $O(a^2_2(a_3-a_2)^2a^3_3)\!=\!O(a^2_2 a^5_3)$, 
and we are done. 
\qed 

We now quote the following important result on lower binomials.  
\begin{thm}
\label{thm:ai}
\cite[Thm.\ 4.5]{ai}
Suppose  $(f,p)\!\in\!\Z[x_1]\times \Pro$, $(v,1)$ is an inner normal
to a lower edge $E$ of $\newt_p(f)$, the lower polynomial $g$ corresponding
to $E$ is a binomial with exponents $\{a_i,a_j\}$, and $p$ does {\em not}
divide $a_i-a_j$. Then the number of roots $\zeta\!\in\!\Q_p$ of
$f$ with $\ord_p \zeta\!=\!v$ is exactly the number of roots of
$g$ in $\Q_p$. \qed
\end{thm}

Finally, we recall a deep theorem from Diophantine approximation 
that allows us to bound from above the $p$-adic valuation of certain high 
degree binomials. 
\begin{yu} 
\scalebox{.95}[1]{\cite[pg.\ 242]{yu} 
Suppose $p\!\in\!\N$ is any prime; $\alpha_1,\ldots, 
\alpha_m$ are nonzero integers;}\linebreak 
\scalebox{.95}[1]{and 
$\beta_1,\ldots,\beta_m$ are integers not all zero. Then  
$\alpha^{\beta_1}_1 \cdots \alpha^{\beta_m}_m\!\neq\!1$ implies that 
$\ord_p\!\left(\alpha^{\beta_1}_1 \cdots \alpha^{\beta_m}_m 
-1\right)$}\linebreak  
\mbox{}\hfill $<22000\left(\frac{9.5(m+1)}{\sqrt{\log p}}\right)^{2(m+1)}
(p-1)\log(10mh) \max\{3,\log\max_i|\beta_i|\} \prod^m_{i=1} 
|\log \alpha_i|$,\hfill \mbox{}\\ 
where $h\!=\!\max\{\log\max_i|\alpha_i|,\log p\}$ and the 
imaginary part of $\log$ lies in $[-\pi,\pi]$. \qed 
\end{yu} 

Let us call any $\newt_p(f)$ such that $f$ has 
no lower $m$-nomials with $m\!\geq\!3$ {\em generic}. 
Oppositely, we call $\newt_p(f)$ {\em flat} if it is a line segment.  
Finally, if $p|(a_i-a_j)$ with $\{a_i,a_j\}$ the exponents of some 
lower binomial of $f$ then we call $\newt_p(f)$ {\em ramified}. 
We will see later that certain ramified cases and flat cases are 
where one begins to see the 
subtleties behind proving $\feas_{\Q_p}(\cF_{1,3})\!\in\!\pp$, 
including the need for Yu's Theorem above. 

\subsection{The Proof of Assertion (1) of Theorem \ref{thm:qp}} 
Our underlying certificate will ultimately be a root  
$\zeta_0\!\in\!\Z/p^\ell\Z$ for $f$ (or a slight variant thereof) 
with $\ell\!=\!O\!\left(p\!\size(f)^8 \right)$. Certain cases 
will actually require such a high power of $p$ and this appears to be 
difficult to avoid. 

Let us write $f(x_1)\!=\!c_1+c_2x^{a_2}_1+c_3 x^{a_3}_1$. 
Just as in Section \ref{sub:proof2}, we may assume $c_1\!\neq\!0$ and 
reduce to certifying roots in $\Z_p$. 
We may also assume that the rightmost (or only) lower edge of $f$ is a 
horizontal line segment at height $0$. (And thus $\ord_pc_1\!\geq\!0$ in 
particular.) This is because we can find the $p$-parts of $c_1,c_2,c_3$ in 
polynomial-time via gcd-free bases (via recursive squaring), compute 
$\newt_p(f)$ in time polynomial in $\size_p(f)$ (via standard convex hull 
algorithms, e.g., \cite{edelsbrunner}), and then rescale $f$ 
without increasing $\size(f)$. More precisely, if $\newt_p(f)$ has no lower 
edges of integral slope then we can immediately conclude that $f$ has no 
roots in $\Q_p$ by Lemma \ref{lemma:newt}. So, replacing $f$ by 
the reciprocal polynomial $f^*$ if necessary, we may assume that the rightmost 
lower edge of $f$ has integral slope and then set $g(x_1)\!:=\!p^{-\ord_pc_2}
f\!\left(p^{\frac{\ord_p(c_2)-\ord_p(c_3)}{a_3-a_2}}x_1\right)$.  
The lower hull of $\newt_p(g)$ then clearly has the desired shape,  
and it is clear that $f$ has a root in $\Q_p$ iff $g$ has a root in $\Q_p$. 
In particular, it is easily checked that $\size(g)\!\leq\!\size(f)$. 

To simplify our proof we will assume that $\gcd(a_2,a_3)\!=\!1$ 
(unless otherwise noted), and recover 
the case $\gcd(a_2,a_3)\!>\!1$ at the very end of our proof. The vanishing 
of $\Delta_\cA(f)$, which can be detected in $\pp$ thanks to 
Corollary \ref{cor:tri}, then determines $2$ cases:  

\smallskip 
\noindent 
\mbox{}\hfill \scalebox{1}[1]{{\bf Case (a): $\pmb{\Delta_\cA(f)\!\neq\!0}$}}
\hfill\mbox{}\\ 
Since $\gcd(a_2,a_3)\!=\!1$ we may clearly assume that   
$p$ divides at most one of $\{a_2,a_3,a_3-a_2\}$. 
The shape of the lower hull of $\newt_p(f)$ (which we've already observed 
can be computed in time polynomial in $\size_p(f)$) then determines $2$ 
subcases:  

\smallskip 
\noindent 
\epsfig{file=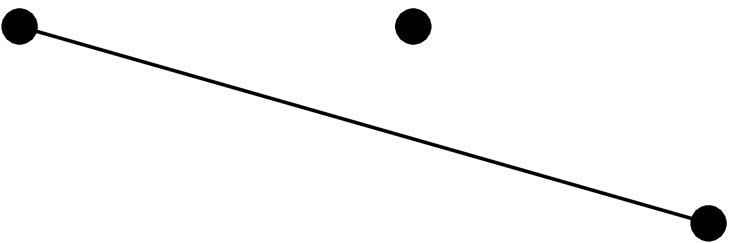,height=.1in} If $\newt_p(f)$ has lower hull a line 
segment then we may also assume (by rescaling $f$ as detailed above) that 
$p\not|c_1,c_3$ and $e\!:=\!\ord_pc_2\!\geq\!0$. 

When $p$ divides either $a_2$ or $a_3-a_2$ 
then we can easily find certificates for  
solvability of $f$ over $\Q_p$: If $e\!=\!0$ then  
$p\not\!|\Delta_\cA(f)$ by Corollary \ref{cor:tri} (since $p\not | a_3$) 
and thus $f$ has no degenerate roots mod $p$. 
So Hensel's Lemma implies that we can use a root of $f$ in $\Z/p\Z$ 
as a certificate for $f$ having a root in $\Q_p$. 
If $e\!>\!0$ then we can in fact detect roots 
in $\Q_p$ for $f$ in $\pp$ by the binomial case, thanks to 
Theorem \ref{thm:ai}. 

So let us now assume 
$p$ does not divide $a_2$ or $a_3-a_2$, and set $e'\!:=\!\ord_pa_3$. 
If $e\!>\!e'$ then observe that 
$f'(x)\!=\!a_3c_3x^{a_3-1}$ mod $p^e$. By Lemma \ref{lemma:newt}, any putative 
root $\zeta\!\in\!\Q_p$ of $f$ must satisfy $\ord_p\zeta \!=\!0$. 
So $f'(\zeta)\!\neq\!0$ mod $p^e$ and Hensel's Lemma implies that 
a root of $f$ in $\Z/p^{2e+1}\Z$ is clearly a 
certificate for $f$ having a root in $\Q_p$. 
Our certificate can also clearly be verified in time polynomial in 
$\size_p(f)$ since $\size(p^{2e+1})\!\leq\!3\size(f)$. 

If $e\!<\!e'$ then 
$f'(x)\!=\!a_2c_2x^{a_2-1}$ mod $p^{e'}$. Similar to  
the last paragraph, $f'(\zeta)\!\neq\!0$ mod $p^{e'}$ and we then 
instead employ a root of $f$ in $\Z/p^\ell\Z$ with 
$\ell\!=\!2e'+1$ as a certificate for $f$ having a root in $\Q_p$. 

Now, if $e\!=\!e'$, observe that 
$\ord_p f'(\zeta)\!=\!\ord_p\frac{f'(\zeta)}{\zeta^{a_2-1}}$ 
since Lemma \ref{lemma:newt} tells us that 
$\ord_p\zeta\!=\!0$\linebreak
\scalebox{.93}[1]{for any root $\zeta\!\in\!\C_p$. Since 
$\Delta_\cA(f)\!\neq\!0$, Corollary \ref{cor:tri} then tells us that 
$\ord_p(a_2c_2+a_3c_3\zeta^{a_3-a_2})\!<\!+\infty$.}\linebreak 
So $\ord_p f'(\zeta)<+\infty$ for any root $\zeta\!\in\!\C_p$ 
of $f$ and then Corollary \ref{cor:tri} tells us that\\ 
\mbox{}\hfill 
$\ord_p \prod_{f(\zeta)=0}f'(\zeta)\!=\!\sum_{f(\zeta)=0}\ord_pf'(\zeta)
\!=\!\ord_p\left((a_3-a_2)^{a_3-a_2}a^{a_2}_2c^{a_3}_2-
(-a_3)^{a_3}c^{a_3-a_2}_1 c^{a_2}_3\right)$\hfill\mbox{}\\
\mbox{}\hfill $=a_3e+\ord_p\left((a_3-a_2)^{a_3-a_2}a^{a_2}_2c^{a_3}_2-
(-a_3)^{a_3}c^{a_3-a_2}_1 c^{a_2}_3\right)$.\hfill \mbox{}\\
(since $p^e|a_2,c_3$).   

So by the $m\!=\!6$ case of Yu's Theorem (using our current assumption that 
$p$ can not divide $a_2$, $a_3-a_2$, $c_1$, or $c_3$),  
we obtain\\
\mbox{}\hfill $\sum_{f(\zeta)=0}\ord_pf'(\zeta)\!=\!a_3e+
O\!\left(p\size(f)^8\right)$. \hfill \mbox{}\\ 
Now, since $p^e|c_2,a_3$, we have $\ord_pf'(\zeta)\!\geq\!e$ for any 
root $\zeta\!\in\!\C_p$ of $f$.  So all roots  
$\zeta\!\in\!\C_p$ of $f$ must satisfy\\ 
\mbox{}\hfill $\ord_pf'(\zeta)\leq e+O\!\left(p\size(f)^8\right)
\leq \size(f)+O\!\left(p\size(f)^8\right)$.\hfill
\mbox{}\\ 
In other words, a root of $f$ in $\Z/p^{O(p\size(f)^8)}\Z$ suffices 
as a certificate, thanks to Hensel's Lemma. 

\medskip 
\noindent 
\epsfig{file=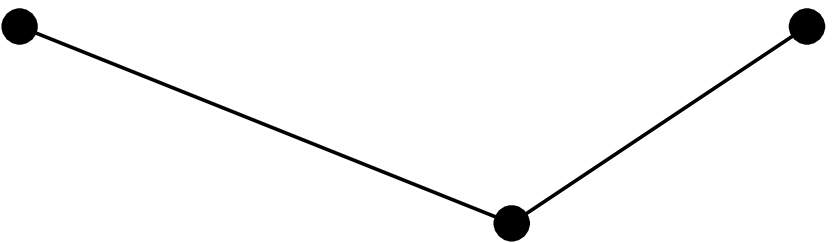,height=.1in} If the lower hull of $\newt_p(f)$ 
is not a line segment then (by rescaling $f$ as detailed above), we may 
also assume that $p|c_1$ but $p\not\!|c_2,c_3$. Since $\gcd(a_2,a_3)\!=\!1$, 
we may also assume (via rescaling and/or reciprocals) that 
$p\not | a_2a_3$, i.e., if $p$ divides the length 
of any lower edge of $\newt_p(f)$ then it is the rightmost 
(now horizontal) edge. 

Via Theorem \ref{thm:ai} and the binomial case of Assertion (1) 
we can easily decide (within time polynomial in $\size_p(f)$) the existence 
of a root of $f$ in $\Z_p$ with valuation $v$, where $(v,1)$ is an 
inner normal of the left lower edge of $\newt_p(f)$. So now we need only 
efficiently detect roots in $\Z_p$ of valuation $0$. Toward this end, 
let us now set $e\!:=\!\ord_p c_1$  and $e'\!:=\!\ord_p(a_3-a_2)$.  
Clearly, $e\!>\!0$ or else we would be in the earlier case where 
$\newt_p(f)$ has lower hull a single edge. 

If $e\!>\!e'$ then $f(x)\!=\!c_2x^{a_2}+c_3x^{a_3}$ mod $p^{e}$ and 
thus $f'(\zeta)\!=\!a_2c_2\zeta^{a_2-1}+a_3c'_3\zeta^{a_3-1}\!=\!-a_2c_3
\zeta^{a_3-1}+a_3c_3\zeta^{a_3-1}\!=\!c_3(a_3-a_2)\zeta^{a_3-1}$  
mod $p^{e}$ for any root $\zeta\!\in\!\C_p$ of $f$. 
So $f'(\zeta)\!\neq\!0$ mod $p^{e}$ for any 
root $\zeta\!\in\!\Z_p$ of valuation $0$ and thus, by  
Hensel's Lemma, we can certify the existence of such a $\zeta$ in 
$\np$ by a root of $f$ in $\Z/p^{2e+1}\Z$. 

If $e\!<\!e'$ then $f'(x)\!=\!a_2c_2x^{a_2-1}+a_3c_3x^{a_3-1}\!=\!a_3c_2
x^{a_2-1}+a_3c_3x^{a_3-1}$ mod $p^{e'}$ since $a_3\!=\!a_2$ mod 
$p^{e'}$. So $f'(\zeta)\!=\!a_3c_2\zeta^{a_2-1}-a_3(c_1\zeta^{-1}
+c_2\zeta^{a_2-1})\!=\!
-\frac{a_3c_1}{\zeta}\!\neq\!0$ mod $p^{e'}$ for any root 
$\zeta\!\in\!\C_p$ of $f$. 
So a root of $f$ in $\Z/p^{2e'+1}\Z$ serves as a certificate for 
a root of $f$ in $\Z_p$. 

Finally, if $e\!=\!e'$, observe that
$f'(x)\!=\!a_2c_2x^{a_2-1}+a_3c_3x^{a_3-1}$ and there are 
exactly $a_2$ (resp.\ $a_3-a_2$) roots of $f$ in $\C_p$ of 
valuation $\frac{e}{a_2}$ (resp.\ $0$) by Lemma \ref{lemma:newt}. 
Using the fact that $p\not | a_2a_3c_2c_3$, it is then easy to 
see that $\ord_p f'(\zeta)=\left(\frac{a_2-1}{a_2}\right)e$ for any 
root $\zeta\!\in\!\C_p$ of $f$ with valuation $\frac{e}{a_2}$. 

The value of $\ord_p f'(\zeta)$ is harder to control at a root of 
valuation $0$. So let us first observe the following:\\ 
($\star$) \hfill 
\scalebox{.9}[1]
{$\frac{a_3c_1}{\zeta}+f'(\zeta)\!=\!\frac{a_3c_1}{\zeta}+a_2c_2\zeta^{a_2-1}
+a_3c_2\zeta^{a_3-1}\!=\!\frac{a_3c_1}{\zeta}+a_3c_2\zeta^{a_2-1}+
a_3c_2\zeta^{a_3-1}\!=\!\frac{a_3}{\zeta}f(\zeta)\!=\!0$ mod $p^e$,}\linebreak 
for any root $\zeta\!\in\!\C_p$ of $f$ of valuation $0$. In other words, 
$e\!\leq\!\ord_p f'(\zeta)$ at any such root. So, similar to 
our earlier flat case, Part (1) of Corollary \ref{cor:tri} implies 
the following:\\
\mbox{}\hfill 
$\ord_p \Delta_\cA(f)\!=\!-(a_2-1)e+\sum\limits_{f(\zeta)=0} f'(\zeta)
\!=\!\sum\limits_{\substack{f(\zeta)=0\\ \ord_\zeta=0}} f'(\zeta)
$.\hfill\mbox{}\\ 
On the other hand, since $e\!=\!\ord_p(a_3-a_2)\!=\!\ord_p c_1$, 
Part (0) of Corollary \ref{cor:tri} combined with the $m\!=\!6$ case of 
Yu's Theorem implies that 
$\ord_p \Delta_\cA(f)\!=\!(a_3-a_2)e+O(p\size(f)^8)$. So any root 
$\zeta\!\in\!\C_p$ of $f$ having valuation $0$ must satisfy\\ 
\mbox{}\hfill 
$\ord_p f'(\zeta)\!\leq\!e+O(p\size(f)^8)\!\leq\!\size(f)+O(p\size(f)^8)$. 
\hfill \mbox{}\\
So again, a root of $f$ in $\Z/p^{O(p\size(f)^8)}\Z$ suffices
as a certificate, thanks to Hensel's Lemma.

\begin{rem} 
Note that if $\newt_p(f)$ is unramified as well as generic, 
then Theorem \ref{thm:ai} implies that we can in fact decide the existence 
of roots in $\Q_p$ for $f$ in $\pp$. \dia 
\end{rem} 

\smallskip
\noindent
\mbox{}\hfill{\bf Case (b): $\pmb{\Delta_\cA(f)\!=\!0}$}\hfill\mbox{}\\
First note that, {\em independent of $\gcd(a_2,a_3)$}, 
a degenerate root of $f$ in $\Q_p$ admits a very 
simple certificate: a $\zeta\!\in\!\Z/p^{4\size(f)+1}\Z$ satisfying  
$c_2(a_3-a_2)\zeta^{a_2}-c_1a_3\!=\!c_3(a_3-a_2)\zeta^{a_3}-c_1a_2
\!=\!0$ mod $p^{4\size(f)+1}$. Thanks to Corollary \ref{cor:tri} and 
our proof of Assertion (0) in Section \ref{sec:binom}, it is clear that the 
preceding $2\times 1$ binomial system has a solution iff  
$f$ has a degenerate root in $\Q_p$. 

So now we resume our assumption that $\gcd(a_2,a_3)\!=\!1$ and build 
certificates for the {\em non-degenerate} roots of $f$ in $\Z_p$. 
Toward this end, observe that the proof of Corollary 
\ref{cor:tri} tells us that the unique degenerate root $\zeta$ of 
$f$ lies in $\Q^*$ and satisfies 
$\left[c_1,c_2\zeta^{a_2},c_3\zeta^{a_3}\right]\!=\! 
\gamma[a_3-a_2,-a_3,a_2]$ for some $\gamma\!\in\!\Q$. Clearly then,   
$q(x_1)\!=\!\frac{1}{\gamma}f(\zeta x_1)$, and $f$ has exactly 
the same number of roots in $\Q_p$ as $q$ does.  

So we can henceforth restrict to the special case 
$c_1\!=\!a_3-a_2$, $c_2\!=\!-a_3$, $c_3\!=\!a_2$, and let 
$r(x_1)\!:=\!\frac{f(x_1)}{(x-1)^2}$ and 
$\Delta\!:=\!\Delta_{\{0,\ldots,a_3-2\}}(r)$. 
Should $p\not\!|a_2a_3(a_3-a_2)$ then $f$ is clearly flat 
and thus all the roots of $f$ have valuation $0$. Part (4) of  
Corollary \ref{cor:tri} 
tells us that $\ord_p \Delta\!\leq\!\log_p\!\left((a_2a_3(a_3-a_2))^{O(1)}
\right)\!=\!O(\log(a_2)+\log(a_3))\!=\!O(\size(f))$ and thus 
the product formula from Lemma \ref{lemma:syl} implies that  
$\ord_p r'(\zeta)\!=\!O(\size(f))$ at any root $\zeta\!\in\!\C_p$ 
of $r$. So a root $\zeta_0\!\in\!\Z/p^{O(\size(f))}\Z$ of $r$  
suffices as a certificate for $f$ to have a root in $\Q_p$ 
other than $1$. (Note also that by construction, $r$ can clearly be evaluated 
mod $p^{O(\size(f))}$ within a number of arithmetic operations 
quadratic in $\size_p(f)$.) 

So let us assume that $p$ divides exactly one number from 
$\{a_2,a_3,a_3-a_2\}$. (Otherwise, $p$ would divide all $3$ numbers, thus 
contradicting the assumption $\gcd(a_2,a_3)\!=\!1$.) 

\smallskip 
\noindent 
\epsfig{file=seglow.eps,height=.1in} 
Should $p|a_3$ then $f$ is clearly flat and, by Lemma \ref{lemma:newt}, every 
root of $r$ has valuation $0$. This implies  
$\ord_p r'(\zeta)\!\geq\!0$ at any root $\zeta\!\in\!\C_p$ 
of $r$. So by Part (4) of
Corollary \ref{cor:tri} and the product formula from Lemma \ref{lemma:syl}, 
we obtain that\\
\mbox{}\hfill $\ord_p \Delta \!=\! (a_3-3)\ord_p(a_2)+
\sum\limits_{r(\zeta)=0} \ord_p r'(\zeta)\!=\!(a_3+O(1))\ord_p(a_2)$. 
\hfill\mbox{}\\
So $\ord_p r'(\zeta)\!=\!O(\ord_p a_2)\!=\!O(\size(f))$ at any root 
$\zeta\!\in\!\C_p$ and we can again use 
a root $\zeta_0\!\in\!\Z/p^{O(\size(f))}\Z$ of $r$  
as a certificate for $f$ to have a root in $\Q_p$ other than $1$. 

\smallskip 
\noindent 
\epsfig{file=bowlow.eps,height=.1in} 
Replacing $f$ by the reciprocal polynomial $f^*$ if need be, we are left with 
the case $p|(a_3-a_2)$. By Lemma \ref{lemma:newt}, $f$ clearly has 
exactly $a_2$ (resp.\ $a_3-a_2$) roots of valuation 
$\frac{\ord_p (a_3-a_2)}{a_2}\!>\!0$ (resp.\ $0$) in $\C_p$. 
Observe that $f'(\zeta)\!=\!a_2a_3\zeta^{a_2-1}(\zeta^{a_3-a_2}-1)$.  

For $\zeta\!\in\!\C_p$ a root of $f$ with valuation 
$\frac{\ord_p (a_3-a_2)}{a_2}$ we then obtain\\
\mbox{}\hfill $\ord_p f'(\zeta)=
\frac{a_2-1}{a_2}\ord_p(a_3-a_2)=O(\size(f))$.\hfill\mbox{}\\ 
In other words, we can simply apply Hensel's Lemma to $f$ and use 
a root of $f$ in\linebreak 
$p^{\ord_p(a_3-a_2)/a_2}(\Z/p^{2\ord_p(a_3-a_2)+1}\Z)$ 
as a certificate for a non-degenerate root of $f$ in $\Q_p$. 

For $\zeta\!\in\!\C_p$ a root of $f$ with valuation $0$ we then 
obtain $\ord_p f'(\zeta)\!\geq\!\ord_p(a_3-a_2)$, thanks to 
identity ($\star$) from the non-degenerate case. Note also that 
$r'(\zeta)\!=\!\frac{f'(\zeta)}{(\zeta-1)^2}-2\frac{f(\zeta)}
{(\zeta-1)^3}\!=\!\frac{f'(\zeta)}{(\zeta-1)^2}$.  
Employing the product formula from 
Lemma \ref{lemma:syl} we then obtain\\  
\mbox{}\hfill $\ord_p \Delta \!=\! 
\left(\sum\limits_{r(\zeta)=0} \ord_p f'(\zeta)\right)-2\ord_p 
\prod\limits_{r(\zeta)=0} (\zeta-1)\!=\!\left(\sum\limits_{r(\zeta)=0} 
\ord_p f'(\zeta)\right)-2\ord_p r(1)$ 
\hfill\mbox{}\\ 
since \mbox{$p\not | a_2$}. 
From our proof of Part (4) of Corollary \ref{cor:tri} it easily follows 
that\linebreak $|r(1)|\!\leq\!a_2a_3(a_3-a_2)$ and thus $\ord_p r(1)\!\leq\!
\log_p(a_2a_3(a_3-a_2))$. So, applying Part (4) of Corollary \ref{cor:tri} 
one last time we obtain \\
\mbox{}\hfill $\sum\limits_{r(\zeta)=0} \ord_p f'(\zeta) \leq 
(a_3+O(1))\ord_p(a_3-a_2)+\log_p(a_2a_3(a_3-a_2))$.\hfill \mbox{}\\
and thus \\ 
\mbox{}\hfill $\sum\limits_{\substack{r(\zeta)=0\\ \ord_p\zeta=0}} 
\ord_p f'(\zeta) \leq 
(a_3-a_2+O(1))\ord_p(a_3-a_2)+\log_p(a_2a_3(a_3-a_2))$.\hfill \mbox{}\\
Since $\ord_p f'(\zeta)\!\geq\!\ord_p(a_3-a_2)$ at a valuation 
$0$ root $\zeta\!\in\!\C_p$ of $f$, and there are exactly 
$a_3-a_2$ such roots, we therefore must have\\ 
\mbox{}\hfill 
$\ord_p f'(\zeta)\!=\!O(1)\ord_p(a_3-a_2)+\log_p(a_2a_3(a_3-a_2))\!=\!
O(\size(f))$.\hfill\mbox{}\\ So we can certify non-degenerate roots 
$\zeta\!\in\!\Q_p$ of $f$ with valuation $0$ by a root 
$\zeta_0\!\in\!\Z/p^{O(\size(f))}\Z$ of $r$ mod $p^{O(\size(f))}$ 
not divisible by $p^{\ord_p (a_3-a_2)/a_2}$. 

{\em Wrapping up the case $\gcd(a_2,a_3)\!>\!1$:}  
From our preceding arguments, we see that we are 
left with certifying the existence of {\em non}-degenerate roots in the case 
$g\!:=\!\gcd(a_2,a_3)\!>\!1$. Fortunately, this is simple: we merely find a   
root non-degenerate root 
$\zeta_0\!\in\!\Z/p^{\ell}\Z$ of $\bar{f}\!:=\!c_1+c_2x^{a_2/g}+c_3x^{a_3/g}$ 
as before (with $\ell$ depending on the case $\bar{f}$ falls into), 
{\em also} satisfying the condition that $x^g-\zeta_0$ has a
root in $\Z/p^{\ell}\Z$. Thanks to Corollary \ref{cor:binom}, we are 
done. \qed 

\section{$\np$-hardness in One Variable: Proving Assertions (4) and (5)} 
\label{sec:primes} 
We will first need to develop two key ingredients: (A) Plaisted's 
beautiful connection between Boolean satisfiability and 
roots of unity, and (B) an algorithm for constructing 
moderately small primes $p$ with $p-1$ having many prime factors. 

\subsection{Roots of Unity and NP-Completeness}  
\label{sub:cyclo} 
Let us define $[n]\!:=\!\{1,\ldots,n\}$. 
Recall that any Boolean expression of one of the following forms:\\ 
\mbox{}\hfill $(\diamondsuit)$ $y_i\vee y_j \vee y_k$, \
$\neg y_i\vee y_j \vee y_k$, \
$\neg y_i\vee \neg y_j \vee y_k$,  \
$\neg y_i\vee \neg y_j \vee \neg y_k$, 
with $i,j,k\!\in\![3n]$,\hfill \mbox{} \\
is a $\sat$ {\em clause}. A {\em satisfying assigment} for an 
arbitrary Boolean formula $B(y_1,\ldots,y_n)$ is  
an assigment of values from $\{0,1\}$ to the variables
$y_1,\ldots,y_{n}$ which makes the equality $B(y_1,\ldots,y_n)\!=\!1$ true. 
Let us now refine slightly Plaisted's elegant reduction from $\sat$ to 
feasibility testing for univariate polynomial systems over the 
complex numbers \cite[Sec.\ 3, pp.\ 127--129]{plaisted}. 
\begin{dfn}  
\label{dfn:plai} 
Letting $P\!:=\!(p_{1},\ldots,p_{n})$ denote any 
strictly increasing sequence of primes, let us inductively 
define a semigroup homomorphism $\cP_P$ --- 
the {\em Plaisted morphism with respect to $P$} --- from certain 
Boolean expressions in the variables $y_1,\ldots,y_n$ to $\Z[x]$, as 
follows:\footnote{
Throughout this paper, for Boolean expressions, we will always identify $0$
with ``{\tt False}'' and  $1$ with ``{\tt True}''.} 
(0) $D_P\!:=\!\prod^n_{i=1}p_{i}$, 
(1) $\cP_P(0)\!:=\!1$, (2) $\cP_P(y_i)\!:=\!x^{D_P/p_{i}}-1$,
(3) $\cP_P(\neg B)\!:=\!(x^{D_P}-1)/\cP_P(B)$, for any 
Boolean expression $B$ for which $\cP_P(B)$ has already been 
defined, 
(4) $\cP_P(B_1\vee B_2)\!:=\!\mathrm{lcm}(\cP_P(B_1),\cP_P(B_2))$, 
for any Boolean expressions $B_1$ and $B_2$ for which $\cP_P(B_1)$ and 
$\cP_P(B_2)$ have already been defined. \dia 
\end{dfn} 
\begin{lemma} 
\label{lemma:plai} 
\cite[Sec.\ 3, pp.\ 127--129]{plaisted} 
Suppose $P\!=\!(p_i)^n_{k=1}$ is an increasing sequence of 
primes with $\log(p_{k})\!=\!O(k^\gamma)$ for some constant $\gamma$. Then, 
for all $n\!\in\!\N$ and any clause $C$ of the form 
$(\diamondsuit)$, we have $\size(\cP_P(C))$ polynomial in $n^\gamma$. 
In particular, 
$\cP_P$ can be evaluated at any such $C$ in time polynomial in $n$. 
Furthermore, if $K$ is any field possessing $D_P$ distinct ${D_P}\thth$ 
roots of unity, then a $\sat$ instance $B(y)\!:=C_1(y)\wedge \cdots \wedge 
C_k(y)$ has a satisfying assignment iff the univariate polynomial system 
$F_B\!:=\!(\cP_P(C_1), \ldots,\cP_P(C_k))$ has a root $\zeta\!\in\!K$ 
satisfying $\zeta^{D_P}-1$. \qed 
\end{lemma} 

\noindent
Plaisted actually proved the special case $K\!=\!\C$ of the above lemma, in 
slightly different language, in \cite{plaisted}. 
However, his proof extends verbatim to the more general 
family of fields detailed above. 

A simple consequence of the resultant is that vanishing
at a $D\thth$ root of unity is algebraically the same thing over $\C$
or $\Q_p$, provided $p$ lies in the right arithmetic progression.

\begin{lemma}
\label{lemma:uni} 
Suppose $D\!\in\!\N$, $f\!\in\!\Z[x]$, and
$p$ is any prime congruent to $1$ mod $D$.
Then $f$ vanishes at a complex $D\thth$ root of unity $\Longleftrightarrow
f$ vanishes at a $D\thth$ root of unity in $\Q_p$. 
\end{lemma} 
\begin{rem} 
Note that $x^2+x+1$ vanishes at a $3\rd$ root of unity in $\C$,
but has {\bf no} roots at all in $\F_5$ or $\Q_5$. So our
congruence assumption on $p$ is necessary. \dia  
\end{rem} 

\smallskip
\noindent
{\bf Proof of Lemma \ref{lemma:uni}:}
First note that by our assumption
on $p$, $\Q_p$ has $D$ distinct $D\thth$ roots of
unity: This follows easily from Hensel's Lemma  
and $\F_p$ having $D$ distinct $D\thth$ roots of unity.
Since $\Z\hookrightarrow\Q_p$ and $\Q_p$
contains all $D\thth$ roots of unity by construction, the
equivalence then follows directly from Lemma 2.8.  \qed

\subsection{Randomization to Avoid Riemann Hypotheses: Proving\\ 
Theorem \ref{thm:primes}}
\label{sub:agp} 
The result below allows us to prove Theorem \ref{thm:primes} and further 
tailor Plaisted's clever reduction to our purposes.  
We let $\pi(x)$ denote the number of primes $\leq\!x$, and 
let $\pi(x;M,1)$ denote the number of primes $\leq\!x$ that are congruent 
to $1 \; \mod \; M$. 

\begin{agp} 
(very special case of \cite[Thm.\ 2.1, pg.\ 712]{carmichael})  
There exist $x_0\!>\!0$ and an $\ell\!\in\!\N$ such that for each  
$x\!\geq\!x_0$, there is a subset $\cD(x)\!\subset\!\N$ of finite cardinality 
$\ell$ with the following property: If $M\!\in\!\N$ satisfies 
$M\!\leq\!x^{2/5}$ and $a\!\not|M$ for all $a\!\in\!\cD(x)$ then 
$\pi(x;M,1)\!\geq\!\frac{\pi(x)}{2\varphi(M)}$. \qed 
\end{agp} 

\noindent 
For those familiar with \cite[Thm.\ 2.1, pg.\ 712]{carmichael}, 
the result above follows immediately upon specializing the  
parameters there as follows:\\ 
\mbox{}\hfill $(A,\eps,\delta,y,a)\!=\!(49/20,1/2,2/245, x,1)$ \hfill 
\mbox{}\\ 
(see also \cite[Fact 4.9]{von}). 

The AGP Theorem enables us to construct random primes from 
certain arithmetic progressions with high probability. An additional 
ingredient that will prove useful is the famous {\em AKS algorithm} 
for deterministic polynomial-time primality checking \cite{aks}. 
Consider now the following algorithm. 
\begin{algor}\mbox{}\\ 
\label{algor:primes}{\bf Input:} A constant $\delta\!>\!0$, a failure 
probability $\eps\!\in\!(0,1/2)$, a positive integer $n$, and the 
constants $x_0$ and $\ell$ from the AGP Theorem.\\ 
{\bf Output:} An increasing sequence $P\!=\!(p_j)^n_{j=1}$ 
of primes, and $c\!\in\!\N$, such that 
$p\!:=\!1+c\prod^n_{i=1} p_i$ satisfies 
$\log p\!=\!O(n\log(n)+\log(1/\eps))$ and, 
with probability $1-\eps$, $p$ is prime. In particular, the output always 
gives a true declaration as to the primality of $p$. 

\medskip 
\noindent
{\bf Description:}\\
0.\ Let $L\!:=\!\lceil 2/\eps\rceil\ell$ 
and compute the first $nL$ primes $p_1, \ldots,p_{nL}$ in increasing order. \\ 
1.\ Define (but do not compute) $M_j\!:=\!\prod^{jn}_{k=(j-1)n+1} p_k$ 
for any $j\!\in\!\N$. Then compute $M_L$, $M_i$\linebreak 
\mbox{}\hspace{.6cm}for a uniformly random $i\!\in\![L]$, 
and $x\!:=\!\max\left\{x_0,
17,
1+M^{5/2}_L
\right\}$. \\ 
2.\ \scalebox{.97}[1]{Compute $K\!:=\!\lfloor (x-1)/M_i\rfloor$ 
and $J\!:=\!\lceil 2\log(2/\eps)\log x\rceil$.}\\ 
3.\ Pick uniformly random $c\!\in\![K]$ until one either has $p:=1+cM_i$ 
prime, or one has $J$\linebreak 
\mbox{}\hspace{.6cm}such numbers that are each composite (using primality 
checks via the AKS algorithm\linebreak
\mbox{}\hspace{.6cm}along the way).\\ 
4.\ If a prime $p$ was found then output\\ 
\mbox{}\hfill ``{\tt $1+c\prod^{in}_{j=(i-1)n+1}p_j$ 
is a prime that works!}''\hfill\mbox{}\\ 
\mbox{}\hspace{.4cm}and stop. Otherwise, stop and output\\
\mbox{}\hspace{.4cm}``{\tt I have failed to find a suitable prime. 
Please forgive me.}'' \dia 
\end{algor}

\begin{rem} 
In our algorithm above, it suffices to find integer approximations to the 
underlying logarithms and square-roots. In particular, we restrict to 
algorithms that can compute the $\log_2 \cL$ most significant bits of 
$\log \cL$, and the $\frac{1}{2}\log_2 \cL$ most significant bits of 
$\sqrt{\cL}$, using $O((\log \cL)(\log \log \cL)\log \log \log \cL)$  
bit operations. Arithmetic-Geometric Mean Iteration and (suitably tailored) 
Newton Iteration are algorithms that respectively satisfy our 
requirements (see, e.g., \cite{dan} for a detailed description). \dia 
\end{rem} 

\smallskip 
\noindent 
{\bf Proof of Theorem \ref{thm:primes}:} 
It clearly suffices to prove that Algorithm \ref{algor:primes} 
is correct, has a success probability that is at least $1-\eps$, and 
works within\\ 
\mbox{}\hfill $O\!\left(\left(\frac{n}{\eps}\right)^{\frac{3}{2}+\delta}+
(n\log(n)+\log(1/\eps))^{7+\delta}\right)$\hfill\mbox{}\\ 
randomized bit operations, 
for any $\delta\!>\!0$. These assertions are proved directly below. \qed  

\smallskip 
\noindent
{\bf Proving Correctness and the Success Probability  
Bound for Algorithm \ref{algor:primes}:} 
First observe that $M_1,\ldots,M_{L}$ are relatively prime. So at most
$\ell$ of the $M_i$ will be divisible by elements of $\cD(x)$.
Note also that $K\!\geq\!1$ and
$1+cM_i\!\leq\!1+KM_i\!\leq\!1+((x-1)/M_i)M_i\!=\!x$
for all $i\!\in\![L]$ and $c\!\in\![K]$.

Since $x\!\geq\!x_0$ and
$x^{2/5}\!\geq\!(x-1)^{2/5}\!\geq\!\left(M^{5/2}_i\right)^{2/5}\!=\!M_i$ for
all $i\!\in\![L]$, the AGP Theorem implies that with probability
at least $1-\frac{\eps}{2}$ (since $i\!\in\![\lceil 2/\eps\rceil \ell]$ is
uniformly random), the arithmetic
progression $\{1+M_i,\ldots,1+KM_i\}$ contains at
least $\frac{\pi(x)}{2\varphi(M_i)}\!\geq\!\frac{\pi(x)}{2M_i}$ primes.
In which case, the proportion of numbers in $\{1+M_i,\ldots,1+KM_i\}$ that
are prime is $\frac{\pi(x)}{2KM_i}\!>\!\frac{\pi(x)}{2+2KM_i}\!>\!
\frac{x/\log x}{2x}\!=\!\frac{1}{2\log x}$, since
$\pi(x)\!>\!x/\log x$ for all $x\!\geq\!17$ \cite[Thm.\ 8.8.1, pg.\ 233]{bs}.
So let us now assume that $i$ is fixed and $M_i$ is not divisible by
any element of $\cD(x)$.

Recalling the inequality $\left(1-\frac{1}{t}\right)^{ct}\!\leq\!e^{-c}$
(valid for all $c\!\geq\!0$ and $t\!\geq\!1$), we then see that
the AGP Theorem implies that the probability of {\em not} finding a prime
of the form $p\!=\!1+cM_i$ after picking $J$ uniformly random
$c\!\in\![K]$ is bounded above by 
$\left(1-\frac{1}{2\log x}\right)^J \!\leq\!\left(1-\frac{1}{2\log x}
\right)^{2\log(2/\eps)\log x}\!\leq\!e^{-\log(2/\eps) }\!=\!
\frac{\eps}{2}$.

In summary, with probability $\geq\!1-\frac{\eps}{2}
-\frac{\eps}{2}\!=\!1-\eps$, Algorithm \ref{algor:primes}
picks an $i$ with $M_i$ not divisible by any element of $\cD(x)$ and
a $c$ such that $p\!:=\!1+cM_i$ is prime. In
particular, we clearly have that\\ 
\mbox{}\hfill 
$\log p \!=\!O(\log(1+KM_i))\!=\!O(n\log(n)+\log(1/\eps))$. \qed
\hfill \mbox{}

\medskip
\noindent
{\bf Complexity Analysis of Algorithm \ref{algor:primes}:} Let $L'\!:=\!nL$
and, for the remainder of our proof, let $p_{i}$ denote the $i\thth$ prime.
Since $L'\!\geq\!6$, we have that\\ 
\mbox{}\hfill $p_{L'}\!\leq L'(\log(L') + \log \log L')$ \hfill \mbox{}\\
by \cite[Thm.\ 8.8.4, pg.\ 233]{bs}. Recall that the
primes in $[\cL]$ can be listed simply by deleting all multiples of
$2$ in $[\cL]$, then deleting all multiples of $3$ in $[\cL]$,
and so on until one reaches multiples of $\lfloor \sqrt{\cL}\rfloor$.
(This is the classic sieve of Eratosthenes.) Recall also that one
can multiply an integer in $[\mu]$ and an integer $[\nu]$
within\\ 
\mbox{}\hfill $O((\log\mu)(\log \log\nu)(\log\log\log \nu)
+(\log\nu)(\log \log\mu) \log\log\log \mu)$\hfill\mbox{}\\  
bit operations (see, e.g., \cite[Table 3.1, pg.\ 43]{bs}). So let us
define the function $\lambda(a):=(\log\log a)\log\log\log a$.

\noindent
{\bf Step 0:} By our preceding observations, it is easily checked that
Step 0 takes $O(L'^{3/2}\log^3 L')$ bit operations.

\noindent
{\bf Step 1:} This step consists of $n-1$ multiplications of
primes with $O(\log L')$ bits (resulting in $M_L$, which has
$O(n\log L')$ bits),
multiplication of a small power of $M_L$ by a square root of $M_L$,
division by an integer with $O(n\log L')$ bits, a constant number of 
additions of integers of
comparable size, and the generation of $O(\log L)$ random bits.
Employing Remark 2.4 along the way, we thus arrive routinely at
an estimate of \\
\mbox{}\hfill
$O\left(n^2(\log L')\lambda(L')+\log(1/\eps)\lambda(1/\eps))
\right)$ \hfill \mbox{}\\
for the total number of bit operations needed for Step 1.

\noindent
{\bf Step 2:} Similar to our analysis of Step 1, we see that
Step 2 has bit complexity\\
\mbox{}\hfill $O((n\log(L')+\log(1/\eps))\lambda(n\log L'))$. \hfill
\mbox{}

\noindent
{\bf Step 3:} This is our most costly step: Here, we require\\
\mbox{}\hfill $O(\log K)\!=\!O(n\log(L')+\log(1/\eps))$\hfill\mbox{}\\
random bits and $J\!=\!O(\log x)\!=\!O(n\log(L')+\log(1/\eps))$ 
primality tests on integers with\\
\mbox{}\hfill$O(\log(1+cM_i))\!=\!O(n\log(L')+\log(1/\eps))$\hfill\mbox{}\\ 
bits. By an improved version of the AKS primality testing algorithm 
\cite{aks,lp} (which takes $O(N^{6+\delta})$ bit operations to test an $N$ 
bit integer for primality), Step 3 can then clearly be done within\\
\mbox{}\hfill
$O\!\left((n\log(L')+\log(1/\eps))^{7+\delta}\right)$
\hfill\mbox{}\\
bit operations, and the generation of $O(n\log(L')+\log(1/\eps))$ random
bits.

\noindent
{\bf Step 4:} This step clearly takes time on the order of the number of
output bits, which is just $O(n\log(n)+\log(1/\eps))$ as already observed
earlier.

\smallskip
\noindent
{\bf Conclusion:} We thus see that Step 0 and
Step 3 dominate the complexity of our algorithm, and we are left with an
overall randomized complexity bound of\\
\mbox{}\hfill
$O\!\left(L'^{3/2}\log^3(L')+ \left(n\log(L')+\log(1/\eps)\right)^{7+\delta}
\right)$\hfill\mbox{}\\
\mbox{}\hfill $=O\!\left(\left(\frac{n}{\eps}\right)^{3/2}\log^3(n/\eps)
+\left(n\log(n)
+\log(1/\eps) \right)^{7+\delta}\right)$\hfill\mbox{}\\
\mbox{}\hfill $=O\!\left(\left(\frac{n}{\eps}\right)^{\frac{3}{2}+\delta}+
\left(n\log(n)+\log(1/\eps)\right)^{7+\delta}
\right)$\hfill\mbox{}\\
randomized bit operations. \qed

\subsection{The Proof of Assertion (4)} 
We will prove a ($\zpp$) randomized polynomial-time 
reduction from $\sat$ to\linebreak 
$\fqp(\Z[x]\times\Pro)$, making use of the intermediate input families
$\{(\Z[x])^k\; | \; k\!\in\!\N\}\times \Pro$ and $\Z[x]\times\{x^D-1\;
| \; D\!\in\!\N\}\times \Pro$ along the way.

Toward this end, suppose $B(y)\!:=\!C_1(y)\wedge\cdots\wedge C_k(y)$
is any $\sat$ instance. The polynomial system
$(\cP_P(C_1),\ldots,\cP_P(C_k))$, for
$P$ the first $n$ primes (employing Lemma \ref{lemma:plai}), 
then
clearly yields $\feas_\C(\{(\Z[x])^k\; | \;
k\!\in\!\N\})\!\in\!\pp \Longrightarrow \pp\!=\!\np$.
Composing this reduction with Proposition 2.6, we then 
immediately obtain $\feas_\C(\Z[x]\times\{x^D-1\; | \;
D\!\in\!\N\})\!\in\!\pp\Longrightarrow \pp\!=\!\np$.

We now need only find a means of transferring from $\C$ to
$\Q_p$. This we do by preceding our reductions above by a judicious
(possibly new) choice of $P$: by applying Theorem \ref{thm:primes}   
with $\eps\!=\!1/3$ (cf.\ Lemma \ref{lemma:uni})
we immediately obtain the implication\\ 
\mbox{}\hfill $\fqp((\Z[x]\times\{x^D-1\; | \;
D\!\in\!\N\})\times \Pro)\!\in\!\zpp\Longrightarrow \np\!\subseteq\!\zpp$.
\hfill \mbox{} 

To conclude, observe that any root $(x,y)\!\in\!\Q^2_p\setminus\{(0,0)\}$
of the quadratic form $x^2-py^2$ must satisfy
$2\ord_p x\!=\!1+2\ord_p y$ (an impossibility). So the only $p$-adic rational
root of $x^2-p y^2$ is $(0,0)$ and we easily
obtain a polynomial-time reduction from\linebreak 
$\fqp((\Z[x]\times\{x^D-1\; | \;
D\!\in\!\N\})\times\Pro)$ to
$\fqp(\Z[x]\times\Pro)$: simply map any\linebreak 
instance $(f(x),x^D-1,p)$ of the former problem to 
$(f(x)^2-(x^D-1)^2p,p)$. So we are done. \qed

\subsection{The Proof of Assertion (5)} 
If we also have the truth of the Wagstaff Conjecture
then we simply repeat our last proof, replacing our
AGP Theorem-based algorithm with a simple brute-force search. More
precisely, letting $D\!:=\!2\cdot 3\cdots p_n$, we simply
test the integers $1+kD$ for primality, starting with $k\!=\!1$ until
one finds a prime. If Wagstaff's Conjecture is true then we
need not proceed any farther than $k\!=\!O\!\left(\frac{\varphi(D)}{D}
\log^2 D\right)$. (Note that $1\!\leq\!\frac{\varphi(D)}{D}\!<\!D$
for all $D\!\geq\!2$.) Using the AKS algorithm,
this brute-force search clearly has (deterministic) complexity 
polynomial in $\log D$ which in turn is polynomial in $n$. \qed

\section{The Final Corollaries} 
\subsection{Proof of Corollary \ref{cor:faster}}
Our proof of Assertion (1) of Theorem \ref{thm:qp} is, in retrospect, a 
polynomial-time reduction from $\fqp(\cF_{1,3})$ to 
$\feas_{\Z/p^\ell\Z}(\cF_{1,3})$ with $\ell\!=\!O(p\size(f)^8)$. 
Combining this reduction with the hypothesis of Corollary \ref{cor:faster} 
then clearly implies that $\feas_{\Q_p}(\cF_{1,3})$ can be solved in time 
polynomial in $p+\size(f)^8$, so we are done. \qed 

\subsection{Proof of Corollary \ref{cor:lots}} 
\label{sub:lots} 
By Lemma \ref{lemma:syl} we know that $\Delta_\cA(f)$ has 
degree at most $2d-1$ in the coefficients of $f$. 
We also know that for any fixed $f\!\in\!T_\cA(H)$, 
$\Delta_\cA(f)$ is an integer as well, and is thus divisible 
by no more than $1+(2d-1)\log(mH))$ primes. (The last assertion follows 
from Lemma \ref{lemma:syl} again, and the elementary fact that 
an integer $N$ has no more than $1+\log N$ distinct prime factors.)  
Recalling that $\pi(x)\!>\!x/\log x$ for all $x\!\geq\!17$ 
\cite[Thm.\ 8.8.1, pg.\ 233]{bs}, we thus obtain that the 
fraction of primes $\leq\!H$ dividing a nonzero $\Delta_\cA(f)$ is 
bounded above by $\frac{1+(2d-1)\log(mH)}{H/\log H}$. 

Now by the Schwartz-Zippel Lemma \cite{schwartz}, $\Delta_\cA(f)$ 
vanishes for at most $(2d-1)m(2H)^{m-1}$ selections of coefficients 
from $\{-H,\ldots,H\}$. In other words, $\Delta_\cA(f)\!=\!0$ for a 
fraction of at most $\frac{(2d-1)m}{2H+1}$ of the polynomials in 
$T_\cA(H)$. 

Combining our last two fractional bounds, we are done. \qed  

\section*{Acknowledgements} 
We thank Jan Denef for pointing out the reference \cite{birchmccann}, 
and Matt Papanikolas and Paula Tretkoff for valuable discussions on 
the Weil Conjectures.  

\bibliographystyle{acm}

\end{document}